\documentclass[12pt]{article}
\usepackage[english]{babel}
\usepackage{amsmath}
\usepackage{amssymb}
\usepackage{amsthm}

        \newtheorem{lemma}{Lemma}[section]
        
        \newtheorem{theorem}[lemma]{Theorem}
        \newtheorem{corollary}[lemma]{Corollary}

        \newtheorem{remark}[lemma]{Remark}

\newcommand{\ov}{\overline}
\newcommand{\bear}{\begin{array}}
\newcommand{\enar}{\end{array}}
\newcommand{\beq}{\begin{equation}}
\newcommand{\eeq}{\end{equation}}
\newcommand{\beqn}{\begin{eqnarray}}
\newcommand{\eeqn}{\end{eqnarray}}
\newcommand{\beit}{\begin{itemize}}
\newcommand{\eeit}{\end{itemize}}

\newcommand{\e}{{\rm e}}

\newcommand{\ds}{\displaystyle}

\newcommand{\ve}{\varepsilon}
\newcommand{\g}{\gamma}
\renewcommand{\l}{\lambda}

\newcommand{\s}{\sigma}
\newcommand{\G}{\Gamma}

\newcommand{\Om}{\Omega}
\renewcommand{\phi}{\varphi}
\renewcommand{\a}{\alpha}
\renewcommand{\b}{\beta}

\newcommand{\til}[1]{\tilde{#1}}
\newcommand{\wtil}[1]{\widetilde{#1}}
\newcommand{\pn}{\par \noindent}
\newcommand{\med}{\medskip}
\newcommand{\qq}{\qquad}
\newcommand{\q}{\quad}
\newcommand{\hookto}{\hookrightarrow}

\newcommand{\media}[1]{\kern 0.4ex {-} \kern -2.0 ex {\int}_{\kern -0.9 ex {#1}}\;}

\newcommand{\PIK}{\textrm{P}(\textrm{H,K})}

\title{Parabolic integrodifferential identification \\
problems related to radial memory kernels I\footnote{Work partially
supported by the Italian Ministero dell'Universit\`a e della
Ricerca Scientifica e Tecnologica (M.U.R.S.T.).}}

\author{Alberto Favaron (Milan), Alfredo Lorenzi (Milan)
\footnote{The authors are members of G.N.A.M.P.A. of the Italian
Istituto Nazionale di Alta Matematica (I.N.d.A.M.)}}

\date{}

\begin{document}
\maketitle \pn

{\bf Abstract.} We are concerned with the problem of recovering the radial
kernel $k$, depending also on time, in the parabolic
integro-differential equation
$$D_{t}u(t,x)={\cal A}u(t,x)\,+\!\int_0^t\!\! k(t-s,|x|)\mathcal{B}u(s,x)ds
\,+\!\int_0^t\!\! D_{|x|}k(t-s,|x|)\mathcal{C}u(s,x)ds+f(t,x),$$
${\cal A}$ being a uniformly elliptic second-order linear
operator in divergence form. We single out a special class of operators 
${\cal A}$
and two pieces of suitable additional information for which the
problem of identifying $k$ can be uniquely solved locally in time
when the domain under consideration is a spherical corona or an annulus.
\med \pn
\med \pn

\section{Posing the identification problem}
\setcounter{equation}{0} The present paper is strictly related to
the previous work \cite{CL} by the latter author and F. Colombo.
Indeed, the problem we are going to
investigate consists in identifying an unknown radial memory
kernel $k$ also depending on time, which appears in the following
integro-differential equation related to the spherical corona
$\Omega\!=\!\{x\!\!=\!\!
(x_1,x_2,x_3)\in\mathbb{R}^3\!:\!R_1<|x|<R_2\}$, where $0<R_1<R_2$
and $|x|={(x_1^2+ x_2^2+ x_3^2)}^{\!1/2}$:
\begin{eqnarray}\label{problem}
D_{t}u(t,x)=\mathcal{A}u(t,x)+\!\int_0^t\!\!
k(t-s,|x|)\mathcal{B}u(s,x)ds+\! \int_0^t\!\!
D_{|x|}k(t-s,|x|)\mathcal{C}u(s,x)ds\;+\!\!\!\!&f(t,x),&
\nonumber\\[2mm]\hskip 8truecm \forall\, (t,x)\in [0,T]
\times\Omega. & &
\end{eqnarray}
In equation $(\ref{problem})\;\mathcal{A}$ and $\mathcal{B}$ are
two second-order linear differential operators, while
$\mathcal{C}$ is a first-order differential operator, having
respectively the following forms:
\begin{eqnarray}
\label{A}
\mathcal{A}=\sum_{j=1}^{3}D_{x_j}\big(\sum_{k=1}^{3}a_{j,k}(x)D_{x_k}
\big),\ \
\mathcal{B}=\sum_{j=1}^{3}D_{x_j}\big(\sum_{k=1}^{3}b_{j,k}(x)D_{x_k}
\big),\ \
\mathcal{C}=\sum_{j=1}^{3}c_{j}(x)D_{x_j}.
\end{eqnarray}
In addition operator $\mathcal{A}$ is uniformly elliptic,\;i.e.
there exist two positive constants ${\alpha}_1$ and ${\alpha}_2$
with ${\alpha}_1\leqslant{\alpha}_2$ such that
\begin{equation}\label{unel}
{\alpha}_1|\xi{|}^2\leqslant\sum_{j,k=1}^{3}a_{j,k}(x){\xi}_j{\xi}_k
\leqslant{\alpha}_2|\xi{|}^2,\qquad\,\forall\,
(x,\xi)\in\Omega\times \mathbb{R}^3\;.
\end{equation}

Before going on, we note that, to the authors' knowledge, the
recover of a kernel $k$ depending also on spatial variables is a
quite new problem, as far as first-order in time
integro-differential equations are concerned. We can quote, besides
\cite{CL}, the papers \cite{JW} and \cite{JJ} that are one-dimensional
in character,
since not only the kernel $k$ is assumed to be {\it degenerate},
i.e. of the form $k(t,x)=\sum_{j=1}^N\,m_j(t)\mu_j(x)$, but also
the space-dependent functions $\mu_j$, $j=1,\ldots,N$, are
assumed to be {\it known}. As a consequence, the identification
problem reduces to recovering the $N$ unknown time-dependent
functions $m_j$, $j=1,\ldots,N $. This latter is nowadays a
classical (vector-) identification problem.

Coming back to our problem, since the  domain $\Omega$ has a
radial symmetry, we will use the classical spherical co-ordinates
$ (r,\varphi,\theta)\in (0,\,\infty)\times (0,\,2\pi) \times
(0,\,\pi)$ related to the Cartesian ones by the well-known
relationship:
\begin{equation}\label{co-ord}
(x_1,x_2,x_3)\!=\!(r\cos{\!\varphi}\sin\!\theta,r\sin{\!\varphi}\sin{\!\theta}
,r\cos{\!\theta})\,.
\end{equation}
Then we prescribe the {\it{initial condition}}
\begin{equation} \label{u0}
u(0,x)=u_0(x)\,,\;\;\;\qquad \forall\, x\in\Omega\,,
\end{equation}
where $u_0\!:\!\overline{\Omega}\rightarrow\mathbb{R}$ is a given
smooth function, as well as one of the following boundary value
conditions, where $u_1\!:\![0,T]\!\times\!{\ov \Om}\to \mathbb{R}$
is a prescribed smooth function, too;
\begin{align}
&\label{DD}(\textrm{D,D})\;\;\quad\qquad u(t,x)=u_1(t,x)\,,
\qquad\qquad\quad\;\;\;\forall\,(t,x)\in
[0,T]\times\partial\mbox{}\Omega\, ,&\\[3mm]
\label{DN} & (\textrm{D,N})\quad\qquad
\left\{\!\!\begin{array}{lll} u(t,x)=u_1(t,x)\,,
\quad\qquad &\forall&\!\!\!\!(t,x)\in
[0,T]\times\partial\mbox{}B(0,R_2)\,,
\\[2,5mm]
\displaystyle\frac{\partial u}{\partial\mbox{}\nu}(t,x)=
\frac{\partial u_1}{\partial\mbox{}\nu}(t,x)\,,\quad\qquad
&\forall&\!\!\!\!(t,x) \in [0,T]\times\partial\mbox{}B(0,R_1)\,,
\end{array}\right.&\\[4mm]
&\label{ND}
(\textrm{N,D})\quad\qquad\left\{\!\!\begin{array}{lll}\displaystyle
\frac{\partial u}{\partial\mbox{}\nu}(t,x)=\frac{\partial
u_1}{\partial\mbox{}\nu} (t,x)\,,\quad\qquad &\forall&\!\!\!\! (t,x)\in
[0,T]\times\partial\mbox{}B(0,R_2)
\,,\\[2,8mm]
\, u(t,x)= u_1(t,x)\,,\qquad\qquad &\forall&\!\!\!\! (t,x)\in
[0,T]\times\partial \mbox{}B(0,R_1)\,,
\end{array}\right.&\\[3mm]
\label{NN} &(\textrm{N,N})\;\;\,\quad\qquad \frac{\partial
u}{\partial\mbox{}\nu}(t,x)=\frac{\partial u_1}{\partial\mbox{}\nu}
(t,x)\,,\qquad\quad\quad\,\forall\, (t,x)\in
[0,T]\times\partial\mbox{}\Omega\,.&
\end{align}
Here D and N stand, respectively, for  the Dirichlet and the conormal
boundary conditions, where the conormal vector $\nu$
is defined by $\nu (x)=\sum_{j,k=1}^3a_{j,k}(x)n_k(x)$, $n(x)$ denoting
the outwarding normal vector at $x\in \partial\mbox{}\Omega$.\\
To determine the radial memory kernel $k$ we need also the two
following additional pieces of information:
\begin{eqnarray}
\label{g1}\!\!\!&\Phi&\!\!\!\!\![u(t,\cdot)](r):=
g_1(t,r)\,,\qquad\forall\,
(t,r)\in[0,T]\times (R_1,R_2)\,,\\[1,3mm]
\label{g2}\!\!\!&\Psi&\!\!\!\!\![u(t,\cdot)]:=
g_2(t)\,,\qquad\qquad\;\forall\, t\in[0,T]\,,
\end{eqnarray}
where $\Phi$ is a linear operator acting on the angular variables
$\varphi,\,\theta$ only, while $\Psi$ is a linear operator acting
on all the
space variables $r,\,\varphi,\,\theta$.\\
\vskip -0,3truecm \pn{\it{Convention:}} from now on we will denote
by $\textrm{P}(\textrm{H,K}),\,\textrm{H,K}\in\{\textrm{D,N}\}$,
the identification problem consisting of $(\ref{problem}),
(\ref{u0})$,
the boundary condition $(\textrm{H,K})$ and $(\ref{g1}),(\ref{g2})$.\\
\vskip -0,3truecm \pn An example of admissible linear operators
$\Phi$ and $\Psi$ is the following:
\begin{align}\label{Phi1}
&\Phi [v](r):= \int_{\!0}^{\pi}\!\!\!\sin\!\theta\mbox{}d\theta\!
\int_{\!0}^{2\pi}\!\!\!\!\lambda(R_2x')v(rx') d\varphi\;,&\\[3mm]
&\label{Psi1}
\Psi[v]:=\int_{\!R_1}^{R_2}\!\!r^2 dr\!\!\int_{\!0}^{\pi}\!\!\!
\sin\!\theta\mbox{} d\theta\!
\int_{\!0}^{2\pi}\!\!\!\!\psi(rx')v(rx')  d\varphi \;,&
\end{align}
where
$x'\!=\!(\cos\!\varphi\sin\!\theta,\,\sin\!\varphi\sin\!\theta,\,
\cos\!\theta) $, while
$\lambda:\partial\mbox{}B(0,R_2)\rightarrow\mathbb{R}$ and
$\psi:\overline{\Omega}\rightarrow\mathbb{R}$ are two smooth
assigned
functions.\\
From $(\ref{DD})\!-\!(\ref{g2})$ we (formally) deduce that our
data have to satisfy the following consistency conditions,
respectively:
\begin{align}\label{DD1}
&(\textrm{C1,D,D})\quad\qquad\quad\q
{u_0}(x)=u_1(0,x),\qquad\qquad\qquad\,\forall\, x\in
\partial\mbox{}\Omega\,,&\\[3mm]
\label{DN1}
&(\textrm{C1,D,N})\qquad\qquad\left\{\!\!\!
\begin{array}{lll} {u_0}(x)=u_1(0,x),
\qquad\qquad\;&\forall&\!\!\!\! x\in \partial\mbox{}B(0,R_2)\;,\\[2,5mm]
\displaystyle\frac{\partial u_0}{\partial\mbox{}\nu}(x)=
\frac{\partial u_1}{\partial\mbox{}\nu}(0,x),\qquad\qquad
&\forall&\!\!\!\! x\in \partial\mbox{}B(0,R_1),\,\,
\end{array}\right.&\\[4mm]
\label{ND1}
&(\textrm{C1,N,D})\qquad\qquad\left\{\!\!\!\begin{array}{lll}
\displaystyle\frac{\partial u_0}{\partial\mbox{}\nu}(x)=
\frac{\partial u_1}{\partial\mbox{}\nu}(0,x),\quad\qquad\quad
&\forall&\!\!\!\! x\in \partial\mbox{}B(0,R_2),\,\\[2,8mm]
\, {u_0}(x)=u_1(0,x),\qquad\qquad\;&\forall&\!\!\!\! x\in
\partial\mbox{}B(0,R_1),\,\,\,
\end{array}\right.&\\[3mm]
\label{NN1}&(\textrm{C1,N,N})\quad\qquad\qquad \frac{\partial
u_0}{\partial\mbox{}\nu}(x)=\frac{\partial u_1}{\partial\mbox{}\nu}(0,x)\,,\qquad\qquad\,\;\forall\,x\in
\partial\mbox{}\Omega\,,&\\\nonumber\\
\label{1.18}&\qq\q\qq\qq\qq\Phi[u_0](r)=g_1(0,r)\,,\quad\quad\,\quad\q\q\forall\,r\in (R_1,R_2)\,,&\\[1,5mm]
\label{1.19}&\qq\q\qq\qq\qq\Psi[u_0]=g_2(0)\,. &\;
\end{align}
Though our identification problem seems to be a very simple
generalization of that in the quoted paper \cite{CL} to the case
where the kernel and the domain $\Om$ are assumed to have radial
symmetries, it should be noted that the choice of $\Om$ coinciding with a ball,
which seems to be the most natural, gives rise to a lot of technical difficulties.
As a consequence, such a problem, in its generality, is still open. Here we stress only
that the mathematical difficulties are concentrated at the centre of the ball
(cf. Remark 2.9).
\pn
However, also in the case of a spherical corona, to solve our
identification problem, we are forced to restrict the {\it admissible}
operators ${\cal A}$ to those satisfying, in addition to $(\ref{A})$,
$(\ref{unel})$
also the following condition for some function $h\in C([R_1,R_2])$:
\begin{eqnarray}\label{RAD}
\sum_{j,k=1}^3\, x_jx_ka_{j,k}(x)=|x|^2h(|x|),\qq \forall x\in
{\ov \Om}.
\end{eqnarray}

\pn
We conclude this section by a remark on radial solutions to our
identification problem.
\begin{remark}
\emph{Assume that ${\cal A}={\cal B}=\Delta_n$, $n=2,3$, ${\cal
C}=D_{|x|}$ and $f$ and $u_0$ are radial functions. We shall refer
to this as to the ``radial case''. If we assumed, as at first
glance seems  reasonable, that in our identification problem also
the state function $u$, i.e. the temperature in physical
applications, should be radial, then definition (\ref{Phi1}) would
reduce to the form
\begin{equation}\label{Phi10}
\Phi [v](r):= C_0v(r),\qq \forall\,r\in [R_1,R_2],
\end{equation}
$C_0$ being a non-zero constant. As a consequence, the additional
condition (\ref{g1}) would amount to requiring that $u$ itself
should be {\it a priori} known. If this is not the case and we
need to determine both $u$ and $k$, we are led to assume that
either of the functions $f$ or $u_0$ is {\it not radial}. \pn On
the contrary, if in the ``radial case'' we {\it a priori} knew a
radial state function $u$, then our problem would reduce to the
following Volterra integro-differential equation of the first kind,
where $n=2,3$:}
\begin{eqnarray}
\label{problemk0} \hskip -1truecm && \int_0^t \big\{ D_r
k(t-s,r)D_r u(s,r) + k(t-s,r)[D_r^2u(s,r) + (n-1)r^{-1}D_r
u(s,r)]\big\}ds = {\til f}(t,r),
\nonumber\\[2mm]
\hskip -1truecm &&\,\hskip 9.5truecm\forall\,
(t,r)\in [0,T]\times [R_1,R_2].
\end{eqnarray}
\emph{Of course the right hand-side $\til{f}$ must satisfy the consistency
condition}
$$
{\til f}(0,r)=0\,,\qq \forall\, r\in [R_1,R_2].
$$
\emph{Furthermore, we note that in this case condition (\ref{g2})
with $\Psi$ being defined by $(\ref{Psi1})$
 makes no
sense, since $\Psi$ reduces to $C_1\int_{R_1}^{R_2}{\til \psi}(r)u(t,r)dr$,
$C_1$ being a constant, i.e. to a {\it known fixed}
function {\it independent} of $k$! \pn However, by
differentiation with respect to time of both sides, equation
(\ref{problemk0}) turns into the equivalent one}
\begin{align}
\label{problemk1}  &\int_0^t \big\{ D_r k(t-s,r)D_tD_r u(s,r) +k(t-s,r)[D_tD_r^2u(s,r)
+ (n-1)r^{-1}D_tD_ru(s,r)]\big\}ds&
\nonumber\\[2mm]
 &\qq + D_r k(t,r)D_r u(0,r) +
k(t,r)[D_r^2u(0,r)+ (n-1)r^{-1}D_ru(0,r)]= D_t{\til f}(t,r)\,,&
\nonumber\\[2,2mm]
 &\hskip 9,4truecm\forall\,
(t,r)\in [0,T]\times [R_1,R_2].&
\end{align}
\emph{To solve this equation we need, e.g., an additional
information of the form}
\begin{eqnarray}
\label{problemk2} \int_{R_1}^{R_2} \l(r)k(t,r)dr=h(t),\qq \forall\,t\in
[0,T].
\end{eqnarray}
\emph{Using the same decomposition for $k$ as in section 3, in section 7 we
will solve the less usual system (\ref{problemk1}) and
(\ref{problemk2}).}
\end{remark}
\section{Main results}
\setcounter{equation}{0} In this section we state our {\it{local
in time}}  existence and uniqueness result related to the
identification problem $\PIK$. For this purpose we assume that the
coefficients of operators
$\mathcal{A},\,\mathcal{B},\,\mathcal{C}\,$ satisfy, in addition
to $(\ref{unel})$ also the following properties:
\begin{align}
\label{ipotesiaij} & a_{i,j}\in W^{2,\infty}({\Omega}),\qquad
a_{i,j}=a_{j,i}\,,\qq\quad i,j=1,2,3\,,&\\
\label{ipotesibijeci} &b_{i,j}\in W^{1,\infty}(\Om)\,,
\qq c_{i}\in L^{\infty}(\Omega)\,,
\quad\quad i,j=1,2,3\,.&
\end{align}
Hence, owing to (1.20) we get that the function
\begin{align}
&\big[\,\widetilde{a}_{1,1}
(r,\varphi,\theta)\,{\cos}^{2}\varphi\,+\,
\widetilde{a}_{1,2}(r,\varphi,\theta)\sin\!2\varphi\,+\,
\widetilde{a}_{2,2}(r,\varphi,\theta)\,{\sin}^{2}\varphi
\big]{\sin}^{2}\theta&\nonumber\\[1mm]
\label{HO}&+\big[\,\widetilde{a}_{1,3}(r,\varphi,\theta)\cos\!\varphi
\,+\,\widetilde{a}_{2,3}(r,\varphi,\theta)\sin\!\varphi
\big]\!\sin\!2\theta\,+\,\widetilde{a}_{3,3}(r,\varphi,\theta)
\,{\cos}^{2}\theta
:= h(r)&
\end{align}
\emph{ depends only on the variable}\, $r$ where we have set
$$
\widetilde{a}_{i,j}(r,\varphi,\theta)\!=\! {a}_{i,j}
(r\cos\!\varphi\sin\!\theta,\,r\sin\!\varphi\sin\!\theta,\,r\cos\!\theta).
$$
\begin{remark}\emph{ To show that our previous condition $(\ref{HO})$ is
meaningful we exhibit a class of coefficients  $a_{i,j}$
satisfying such a property. Let us suppose that there exist three
functions $a,b,d\in W^{2,\infty}(R_1,R_2)$ and a function
$c\in W^{2,\infty}({\Om})$, $a$ and $c$ being, respectively, {\it positive}
and {\it non-negative}, such that}
\begin{equation}\label{condsuaij}
\left\{\!\!\!\begin{array}{lll}
 a_{1,1}(x)\!\!\!&=&\!\!\!a(|x|)+
\displaystyle\frac{(x_2^2+x_3^2)[c(x)-b(|x|)]}{|x|^2}+
\frac{x_1^2d(|x|)}{|x|^2},\\[5mm]
 a_{2,2}(x)\!\!\!&=&\!\!\!a(|x|)+
\displaystyle\frac{(x_1^2+x_3^2)[c(x)-b(|x|)]}{|x|^2}+
\frac{x_2^2d(|x|)}{|x|^2},
\\[5mm]
 a_{3,3}(x)\!\!\!&=&\!\!\!a(|x|)+
\displaystyle\frac{(x_1^2+x_2^2)[c(x)-b(|x|)]}{|x|^2}
+\frac{x_3^2d(|x|)}{|x|^2},
\\[5mm]
 a_{j,k}(x)\!\!\!&= &\!\!\!a_{k,j}(x)=\displaystyle
\frac{x_jx_k[b(|x|)-c(x)+d(|x|)]}{|x|^2},\qq
1\le j,k\le 3,\ j\neq k.
\end{array}\right.
\end{equation}
\emph{Simple computations show that property (\ref{RAD}) is
satisfied with
$h(r)=a(r)+d(r)$.\\
Denoting by $f^+,\, f^-$, respectively, the positive and the
negative parts of  a function $f$, for every $x\in\Omega$ and
$\xi\in\mathbb{R}^3$ from $(\ref{condsuaij})$ it follows}
\begin{align}
& \sum_{j,k=1}^{3}a_{j,k}(x){\xi}_j{\xi}_k=
a(|x|){|\xi|}^2+\frac{c(x)-b(|x|)}{|x|^2}\big[(x_2^2\!
+\!x_3^2){\xi_1}^2+(x_1^2\!+\!x_3^2){\xi_2}^2+(x_1^2+x_2^2){\xi_3}^2&\nonumber\\
&\q\;\;\,- 2x_1x_2\xi_1\xi_2-2x_1x_3\xi_1\xi_3-2x_2x_3\xi_2
\xi_3\big]+\frac{d(|x|)}{|x|^2}{(x_1\xi_1\!+\!x_2\xi_2
+\!x_3\xi_3)}^2\;&
\nonumber\\[1,3mm]
\label{unel1}&\q\;\;\,\geqslant \,a(|x|){|\xi|}^2-\frac{b^+(|x|)}{|x|^2}
{|x\wedge\xi |}^2-\frac{d^-(|x|)}{|x|^2}{[x\cdot\xi]}^2
\ge\,\big[a(|x|)-b^+(|x|)-d^-(|x|)\big]|\xi|^2,&
\end{align}
\emph{where $\wedge$ and $\cdot$ denote, respectively,
the wedge and the inner product in $\mathbb{R}^3$.\\
Hence, to ensure the uniform ellipticity of operator
$\mathcal{A}$, we need the additional assumption:}
\begin{equation}\label{M11}
a(r)-b^+(r)-d^-(r)>0,\qquad\forall\,r\in [R_1,R_2].
\end{equation}
\emph{Condition $(\ref{M11})$ amounts to requiring that function $a$ is
large enough with respect to $b$ and $d^-$.
Consequently $(\ref{unel})$ is trivially
satisfied owing to $(\ref{M11})$ with
$$\a_1\!=\!\min_{r\in [R_1,R_2]}\big[a(r)-b^+(r)-d^-(r)\big]\,,\qq
\a_2\!=\!\|a+b^-+d^+\|_{C([R_1,R_2])}+\|c\|_{C(\ov\Om)}.$$}
\end{remark}
\begin{remark} \emph{We can widen the class of special operators in the previous
remark introducing the following function sequence
$\{a_{i,j}^{(n)}\}_{n=1}^{+\infty}$:
\begin{alignat}{3}
\label{newcond} & a_{1,1}^{(n)}(x)=2x_1^{2n-2}x_2^{2n}x_3^{2n},\q
 a_{2,2}^{(n)}(x)=2x_1^{2n}x_2^{2n-2}x_3^{2n},\q
 a_{3,3}^{(n)}(x)=2x_1^{2n}x_2^{2n}x_3^{2n-2},&\\[2,5mm]
&a_{1,2}^{(n)}(x)=a_{2,1}^{(n)}(x)=-2x_1^{2n-1}x_2^{2n-1}x_3^{2n},&\\[2,5mm]
&a_{1,3}^{(n)}(x)=a_{3,1}^{(n)}(x)=-2x_1^{2n-1}x_2^{2n}x_3^{2n-1},&\\[2,5mm]
&a_{2,3}^{(n)}(x)=a_{3,2}^{(n)}(x)=-2x_1^{2n}x_2^{2n-1}x_3^{2n-1}.&
\end{alignat}
Simple computations show that
\begin{eqnarray}
\label{unel2} \hskip -1.5truecm
\sum_{j,k=1}^{3}\,a_{j,k}^{(n)}(x){\xi}_j{\xi}_k \!\!\!&=&\!\!\!
(x_1x_2x_3)^{2n-2}\big[ (\xi_1x_2x_3\!-\!x_1\xi_2x_3)^2
+(x_1\xi_2x_3\!-\!x_1x_2\xi_3)^2
\nonumber\\
\hskip -1.5truecm && +(x_1x_2\xi_3-\xi_1x_2x_3)^2\big] \ge 0,\qq
\forall \,(x,\xi)\in {\ov \Om}\times \mathbb{R}^3,\ \forall \,n\in
\mathbb{N}.
\end{eqnarray}
In particular we get
\begin{eqnarray}
\label{symm1} && \sum_{j,k=1}^{3}\,x_jx_ka_{j,k}^{(n)}(x)=0,\qq
\forall\, x\in {\ov \Om}, \ \forall\, n\in \mathbb{N}.
\end{eqnarray}
Then with any sequence of {\it non-negative} functions
$\{d_n\}_{n=1}^{+\infty}\subset W^{2,\infty}(\Om)$ we associate the special
coefficients:
\begin{eqnarray}
\label{newcoef} a_{j,k}(x)={\bar a}_{j,k}(x) +
\sum_{n=1}^N\,a_{j,k}^{(n)}(x)d_n(x), \qq x\in {\ov \Om},\
j,k=1,2,3,\ N\in \mathbb{N},
\end{eqnarray}
where the coefficients ${\bar a}_{j,k}$ are defined by equations
(\ref{condsuaij}).\! Also the value $N\!=\!+\infty$ is allowed provided
that the series $\sum_{n=1}^{+\infty}a_{j,k}^{(n)}(x)d_n(x)$
may be differentiated twice term by term, with a sum in
$ W^{2,\infty}(\Om)$.\pn It is immediate to check that both property
(\ref{RAD}) and the uniform ellipticity are satisfied under the
same assumptions as in the previous remark. \pn We have thus
showed that the class of admissible coefficients is not limited to
those represented by (\ref{condsuaij}).}
\end{remark}
\begin{remark}\label{conormal}
\emph{Simple computations show that, when the coefficients $a_{j,k}$
are given by $(\ref{condsuaij})$, then the conormal vector $\nu^{(l)}$
on $\partial\mbox{}B(0,R_l)$ coincides with the normal vector
$$\nu^{(l)}(x)=(-1)^lR_l^{-1}[a(R_l)-d(R_l)]x\,,\qq l=1,2.$$}
\end{remark}
In order to find out the right hypotheses on the linear operators
$\Phi$ and $\Psi$, it will be convenient to rewrite the operator
$\mathcal{A}$ in the spherical co-ordinates
$(r,\,\varphi,\,\theta)$. Recall first that the
gradient $\nabla=(D_{x_1},D_{x_2},D_{x_3})$ can be expressed in
terms of $(D_r,D_{\varphi},D_{\theta})$ by the formulae:
\begin{equation}
\label{D123} \left\{\!\! \begin{array}{lll} D_{x_1}\!\!\! & =
&\!\!\!{\cos\!\varphi\sin\!\theta} D_{r}-\displaystyle
\frac{\sin\!\varphi}{r\sin\!\theta}D_{\varphi}+\displaystyle
\frac{\cos\!\varphi\cos\!\theta}{r}D_{\theta}\,,\\[3,0mm]
D_{x_2}\!\!\!  & = &\!\!\!{\sin\!\varphi\sin\!\theta}
D_{r}+\displaystyle
\frac{\cos\!\varphi}{r\sin\!\theta}D_{\varphi}+\displaystyle
\frac{\sin\!\varphi\cos\!\theta}{r}D_{\theta}\,,\\[3,0mm]
D_{x_3} \!\!\! & = &\!\!\!{\cos\!\theta}  D_{r}-\displaystyle
\frac{\sin\!\theta}{r}D_{\theta}\,.
\end{array}\right.
\end{equation}
As a consequence, simple computations easily yield
\begin{eqnarray}
\label{D1}\sum_{k=1}^{3}a_{1,k}(x)D_{x_k}\!\!\!& = &\!\!\!
f_1(r,\varphi,\theta)D_{r}+\frac{f_2(r,\varphi,\theta)}{r\sin\!\theta}
D_{\varphi}+\frac{ f_3(r,\varphi,\theta)}{r}D_{\theta}\,,\\
\label{D2} \sum_{k=1}^{3} a_{2,k}(x)D_{x_k}\!\!\!& = &\!\!\!
g_1(r,\varphi,\theta)D_{r}+\frac{g_2(r,\varphi,\theta)}{r\sin\!\theta}
D_{\varphi}+\frac{ g_3(r,\varphi,\theta)}{r}D_{\theta}\,,\\
\label{D3} \sum_{k=1}^{3}a_{3,k}(x)D_{x_k}\!\!\!& = &\!\!\!
h_1(r,\varphi,\theta)D_{r}+
\frac{h_2(r,\varphi,\theta)}{r\sin\!\theta}D_{\varphi}+ \frac{
h_3(r,\varphi,\theta)}{r}D_{\theta}\,,\qquad\quad
\end{eqnarray}
functions $f_j,\,g_j,\,h_j,\, j\!=\!1,2,3,$ being defined by
\begin{eqnarray}\label{fj}
\left\{\!\!
\begin{array}{lll}
f_1(r,\varphi,\theta)\!\!\!\!&:=&\!\!\!
{\widetilde{a}}_{1,1}(r,\varphi,\theta)\!
\cos\!\varphi\sin\!\theta+{\widetilde{a}}_{1,2}(r,\varphi,\theta)\!\sin\!\varphi
\sin\!\theta+\widetilde{a}_{1,3}(r,\varphi,\theta)\!\cos\!\theta\,,\\[3mm]
f_2(r,\varphi,\theta)\!\!\!\!&:=&\!\!\!\widetilde{a}_{1,2}(r,\varphi,\theta)\!
\cos\!\varphi-
{\widetilde{a}}_{1,1}(r,\varphi,\theta)\!\sin\!\varphi\,,\qquad
\qquad\qquad\qquad\qquad\quad\;\;\, \\[3mm]
f_3(r,\varphi,\theta)\!\!\!\!&:=&\!\!\!
\widetilde{a}_{1,1}(r,\varphi,\theta)\! \cos\!\varphi\cos\!\theta
+{\widetilde{a}}_{1,2}(r,\varphi,\theta)\!\sin\!\varphi
\cos\!\theta-{\widetilde{a}}_{1,3}(r,\varphi,\theta)\!\sin\!\theta\,,
\end{array}
\right.\\[3mm]
\label{gj} \left\{\!\! \begin{array}{lll}
g_1(r,\varphi,\theta)\!\!\!\!&:=&\!\!\!
{\widetilde{a}}_{2,1}(r,\varphi,\theta)\!
\cos\!\varphi\sin\!\theta+{\widetilde{a}}_{2,2}(r,\varphi,\theta)\!\sin\!\varphi
\sin\!\theta+\widetilde{a}_{2,3}(r,\varphi,\theta)\!\cos\!\theta\,,\\[3mm]
g_2(r,\varphi,\theta)\!\!\!\!&:=&\!\!\!\widetilde{a}_{2,2}(r,\varphi,\theta)\!
\cos\!\varphi-
{\widetilde{a}}_{2,1}(r,\varphi,\theta)\!\sin\!\varphi\,,\qquad\qquad
\qquad\qquad\qquad\quad\;\;\, \\[3mm]
g_3(r,\varphi,\theta)\!\!\!\!&:=&\!\!\!
\widetilde{a}_{2,1}(r,\varphi,\theta)\! \cos\!\varphi\cos\!\theta
+{\widetilde{a}}_{2,2}(r,\varphi,\theta)\!\sin\!\varphi
\cos\!\theta-{\widetilde{a}}_{2,3}(r,\varphi,\theta)\!\sin\!\theta\,,
\end{array}\right.\\[3mm]
\label{hj} \left\{\!\! \begin{array}{lll}
h_1(r,\varphi,\theta)\!\!\!\!&:=&\!\!\!
{\widetilde{a}}_{3,1}(r,\varphi,\theta)\!
\cos\!\varphi\sin\!\theta+{\widetilde{a}}_{3,2}(r,\varphi,\theta)
\!\sin\!\varphi
\sin\!\theta+\widetilde{a}_{3,3}(r,\varphi,\theta)\!\cos\!\theta\,,\\[3mm]
h_2(r,\varphi,\theta)
\!\!\!\!&:=&\!\!\!\widetilde{a}_{3,2}(r,\varphi,\theta)\!
\cos\!\varphi-
{\widetilde{a}}_{3,1}(r,\varphi,\theta)\!\sin\!\varphi\,,\qquad\qquad
\qquad\qquad\qquad\quad\;\;\,  \\[3mm]
h_3(r,\varphi,\theta)\!\!\!\!&:=&\!\!\!
\widetilde{a}_{3,1}(r,\varphi,\theta)\! \cos\!\varphi\cos\!\theta
+{\widetilde{a}}_{3,2}(r,\varphi,\theta)\!\sin\!\varphi
\cos\!\theta-{\widetilde{a}}_{3,3}(r,\varphi,\theta)\!\sin\!\theta.
\end{array}\right.
\end{eqnarray}
Hence, using again $(\ref{D123})$ and applying it to relations
$(\ref{D1})\! -\!(\ref{D3})$, we get:
\begin{align}
&D_{x_1}\Big(\sum_{k=1}^{3}a_{1,k}(x)D_{x_k}
\Big)= D_{r}\Big[ f_1(r,\varphi,\theta)\!\cos\!\varphi\sin\!\theta
D_{r} + \displaystyle\frac{
f_2(r,\varphi,\theta)\!\cos\!\varphi}{r}D_{\varphi}&
\nonumber \\
&\q\q+\displaystyle\frac{
f_3(r,\varphi,\theta)\!\cos\!\varphi
\sin\!\theta}{r}D_{\theta}\Big] -\frac{
\sin\!\varphi}{r\sin\!\theta } D_{\varphi}\Big[
f_1(r,\varphi,\theta)D_r+\frac{ f_2(r,\varphi,\theta)}
{r\sin\!\theta}D_{\varphi}+\frac{
f_3(r,\varphi,\theta)}{r}D_{\theta}\Big]&\nonumber\\[2mm]
\label{D1D1}&\q\q+ \displaystyle\frac{ \cos\!\theta}{r}
D_{\theta} \Big[f_1(r,\varphi,\theta)\!\cos\!\varphi D_{r}
 +\frac{f_2(r,\varphi,\theta)\!\cos\!\varphi}{r\sin\!\theta}D_{\varphi}+
\frac{f_3(r,\varphi,\theta)\!\cos\!\varphi}{r}D_{\theta}\Big]\,,&
\end{align}
\begin{align}
 &D_{x_2}\Big(\sum_{k=1}^{3}
a_{2,k}(x)D_{x_k}  \Big)= D_{r}\Big[
g_1(r,\varphi,\theta)\!\sin\!\varphi\sin\!\theta D_{r} +
\frac{g_2(r,\varphi,\theta)\!\sin\!\varphi}{r}D_{\varphi}&
\nonumber \\
&\q\q+\frac{g_3(r,\varphi,\theta)\!\sin\!\varphi\sin\!\theta}{r}
D_{\theta}\Big] +\frac{ \cos\!\varphi}{r\sin\!\theta }
D_{\varphi}\Big[ g_1(r,\varphi,\theta)D_r+\frac{
g_2(r,\varphi,\theta)}{r\sin\!\theta}
D_{\varphi}+\frac{ g_3(r,\varphi,\theta)}{r}D_{\theta}\Big]&\nonumber\\[2mm]
&\q\q+\frac{ \cos\!\theta}{r}
D_{\theta}\Big[g_1(r,\varphi,\theta)\! \sin\!\varphi D_{r}
 +\frac{ g_2(r,\varphi,\theta)\!\sin\!\varphi}{r\sin\!\theta}D_{\varphi}+
\frac{ g_3(r,\varphi,\theta)\!\sin\!\varphi
}{r}D_{\theta}\Big]\,,&\label{D2D2}
\end{align}
\begin{align}
&D_{x_3}\Big(\sum_{k=1}^{3}  a_{3,k}D_{x_k}
\Big)= D_{r}\Big[ h_1(r,\varphi,\theta)\!\cos\!\theta D_{r}
+\frac{ h_2(r,\varphi,
\theta)\!\cos\!\theta}{r\sin\!\theta}D_{\varphi}+\frac{h_3(r,\varphi,\theta)\!
\cos\!\theta }{r}D_{\theta}\Big]&\nonumber \\
\label{D3D3}&\qq\qq-\frac{ \sin\!\theta}{r} D_{\theta}
\Big[\frac{h_1(r,\varphi,\theta)}{r}D_{r}
 +\frac{ h_2(r,\varphi,\theta)}{r\sin\!\theta}D_{\varphi}
+\frac{h_3(r,\varphi,\theta)}{r}D_{\theta}\Big]\;.&
\end{align}
\vskip 0,1truecm\pn
 Let us now define the following functions, where $j=1,2,3$:
\begin{equation}\label{kj,lj}
\left\{\!\!\!\begin{array}{lll} k_j(r,\varphi,\theta) \!\!\!\!\!
&:\, = & \!\!\!\!f_j(r,\varphi,\theta)\!
\cos\!\varphi\sin\!\theta+
g_j(r,\varphi,\theta)\!\sin\!\varphi\sin\!\theta
+ h_j(r,\varphi,\theta)\!\cos\!\theta\,, \\[2mm]
l_j(r,\varphi,\theta) \!\!\!\!\!  & :\,  = & \!\!\!\!
f_j(r,\varphi,\theta) \!\cos\!\varphi
+g_j(r,\varphi,\theta)\!\sin\!\varphi\,.
\end{array}\right.
\end{equation}
Using $(\ref{fj})\!-\!(\ref{hj})$ one can easily check that $k_1$
coincides with the function $h$ defined in $(\ref{HO})$. Therefore,
by virtue of the assumptions (\ref{RAD}) made on the coefficients
$a_{i,j}$, we conclude that $k_1$
depends on $r$  only.\\
Then, rearranging the terms on the right-hand sides of
$(\ref{D1D1})\! -\!(\ref{D3D3})$, we  obtain the following
polar representation $\widetilde{\mathcal{A}}$ for the second-order
differential operator $\mathcal{A}$:
\begin{eqnarray}
{\widetilde{\mathcal{A}}}\!\!\!&=&\!\!\!
D_r\Big[k_1(r)D_r+\frac{k_2(r,\varphi,\theta)}{r\sin\!\theta}D_{\varphi}+
\frac{k_3(r,\varphi,\theta)}{r}D_{\theta} \Big]- \frac{
\sin\!\varphi}{r\sin\!\theta } D_{\varphi}\Big[
{f_1(r,\varphi,\theta)}D_r\nonumber \\[2mm]
&& +\,\frac{f_2(r,\varphi,\theta)}{r\sin\!\theta}D_{\varphi}  +\,\frac{f_3(r,\varphi,\theta)}{r}D_{\theta}\Big]+
\frac{ \cos\!\varphi} {r\sin\!\theta }
D_{\varphi}\Big[{g_1(r,\varphi,\theta)}
D_r+\frac{g_2(r,\varphi,\theta)}{r\sin\!\theta}D_{\varphi}
\nonumber\\[2mm]
&& +\,\frac{g_3(r,\varphi,\theta)}{r}D_{\theta}\Big] -\, \frac{ \sin\!\theta}{r}
D_{\theta}\Big[{h_1(r,\varphi,\theta)}
D_r+\frac{h_2(r,\varphi,\theta)}{r\sin\!\theta}D_{\varphi}+
\frac{h_3(r,\varphi,\theta)}{r}D_{\theta}\Big]
\nonumber\\[2mm]
\label{tildeA} && +\, \frac{ \cos\!\theta}{r} D_{\theta}
\Big[{l_1(r,\varphi,\theta)}D_r\,+\,\frac{l_2(r,\varphi,\theta)}
{r\sin\!\theta}D_{\varphi}\,+\,\frac{l_3(r,\varphi,\theta)}{r}D_{\theta}\Big]\,.
\end{eqnarray}
We can now list our requirements on  operators $\Phi$ and $\Psi$
in accordance with the explicit case $(\ref{Phi1})$ and
$(\ref{Psi1})$. We will work in Sobolev spaces related to $L^p(\Om)$ with
\begin{equation}\label{p3}
p\in (3,+\infty)
\end{equation}
and we will assume
\begin{alignat}{5}
\label{primasuPhiePsi} &
\Phi\in\mathcal{L}\big(L^p(\Omega),\,L^p(R_1,R_2)\big),\q\q\qq\,\,\,\,
 \Psi\in L^p(\Omega)^*,&\\[1,9mm]
\label{secondasuPhi} & \Phi[wu]=w\,\Phi[u],\qq\qq\qq\qq\q\,\,\,
 \forall\,(w,u)\in L^p(R_1,R_2)\times L^p(\Omega),&\\[1,9mm]
\label{terzasuPhi} & D_r\Phi[u](r)=\Phi[D_ru](r), \qq\qq\q\;\;\;\,\,
\forall\,u\in W^{1,p}(\Omega)\;\;\textrm{and}\,\,r\in(R_1,\,R_2),&
\\[1,7mm]
\label{quartasuPhi} &
\Phi\mathcal{\widetilde{A}}=\mathcal{\widetilde{A}}_1\Phi+{\Phi}_1
\quad\;\textrm{on}\;W^{2,p}(\Omega),\qq\;\,\,\,
 \Phi_1\in \mathcal{L}\big( W^{1,p}(\Omega),L^p(R_1,R_2)\big),&\\[1,9mm]
\label{primasuPsi} &
\Psi\mathcal{\widetilde{A}}={\Psi}_1\quad\quad\textrm{on}\;
W^{2,p}(\Omega),\q\q\q\q\q\,  {\Psi}_1\in W^{1,p}({\Omega})^*,&
\end{alignat}
where
\begin{equation}\label{tildeA1}
\mathcal{\widetilde{A}}_1= D_r[k_1(r)D_r]\;.
\end{equation}
To state our result concerning the identification problem
$(\ref{problem}),\, (\ref{u0})\!-\!(\ref{g2})$ we need to list
also the following  assumptions on the data
$f,\,u_0,\,u_1,\,g_1,\,g_2$:
\begin{align}
\label{richiestasuf}
&f\in C^{1+\beta}\big([0,T];L^p(\Omega)\big)\,,\q\;
f(0,\cdot)\in W^{2,p}(\Omega)&\\[2,3mm]
\label{richiestaperu0}&u_0\in W^{4,p}(\Omega)\;,\q
{\cal B}u_0\in W_{\textrm{H,K}}^{2\delta,p}(\Omega)\;\, &\\[2,3mm]
\label{richiestaperu1} & u_1\in
C^{2+\beta}\big([0,T];L^{p}(\Omega)\big)\cap
C^{1+\beta}\big([0,T];W^{2,p}(\Omega)\big)\,,&\\[2,3mm]
\label{richiestaperAu0}&\mathcal{A}u_0+f(0,\cdot)-
D_tu_1(0,\cdot)\in W_{\!\textrm{H,K}}^{2,p}(\Omega)\,,&\\[2,3mm]
\label{richiestaperA2u0}&F\!:=k_0'{\cal C}u_0+k_0{\cal B}u_0 +{\cal A}^2u_0+{\cal A}f(0,\cdot)
-D_t^2u_1(0,\cdot)+D_t f(0,\cdot)\in
W_{\textrm{H,K}}^{2\b,p}(\Omega)\,,&\\[2,3mm]
\label{richiesteperg1} &g_1\in
C^{2+\beta}\big([0,T];L^{p}(R_1,R_2)\big)\cap C^{1+\b}\big([0,T];
W^{2,p}(R_1,R_2)\big)\,,&\\[2,3mm]
\label{richiesteperg2}
&g_2 \in
C^{2+\beta}\big([0,T];\mathbb{R}\big)\,,&
\end{align}
where $\beta\!\in \!(0,1/2)\backslash \{1/(2p)\}$, $\delta \!\in\!(\b,1/2)\backslash \{1/(2p)\}$,
and function $k_0$ in $(\ref{richiestaperA2u0})$ is defined by formula $(3.18)$.
Moreover, the spaces $W_{\!\textrm{H,K}}^{2,p}(\Omega)$, $\textrm{H,K}\in \{\textrm{D,N}\}$,
 are defined by
\begin{equation}\label{WIK}
 W_{\!\textrm{H,K}}^{2,p}(\Omega)\!=\!\big\{\! w \in  W^{2,p}(\Omega)\!:w \;
\textrm{satisfies the homogeneous condition (H,K)}\!\},
\end{equation}
whereas the spaces
$W_{\textrm{H,K}}^{2\g,p}(\Omega)\!\!\equiv
\!\!{\big(L^p(\Omega),
W_{\!\textrm{H,K}}^{2,p}(\Omega)\big)}_{\g,\infty}$,
$\g\!\in \!(0,1/2]\backslash \{1/(2p)\}$,  are
interpolation spaces between $W_{\!\textrm{H,K}}^{2,p}(\Omega)$
and $L^p(\Omega)$ and they are defined (cf. \cite[section 4.3.3]{TR}),
respectively, by the following equations:
\begin{align}
\hskip -5truemm
\label{WDD}&W_{\textrm{D,D}}^{2\g,p}(\Omega)=
\left\{\!\!\begin{array}{lll} W^{2\g,p}(\Omega), & &\qq\;\;\;\;
\textrm{if}\;\; 0<\g<{1}/{(2p)}\,,
\\[1,5mm]
\big\{u\in W^{2\g,p}(\Omega):u=0
\;\;\textrm{on}\;\partial\mbox{}\Omega \big\},&
&\qq\;\;\;\;\textrm{if}\;\;{1}/{(2p)}<\g\le 1/2\,,
\end{array}\right.\\[4,5mm]
\label{WDN}&W_{\textrm{D,N}}^{2\g,p}(\Omega)=
\left\{\!\!\begin{array}{l}
W^{2\g,p}(\Om),\q\qq\qq\qq\qq\qq\qq\q\,
\textrm{if}\;\; 0<\g<{1}/{(2p)}\,,\\[1,6mm]
\big\{u\in W^{2\g,p}(\Omega):u=0 \;\;\textrm{on}\;
\partial\mbox{}B(0,R_2)
\big\},\q\;\;\,\textrm{if}\;\;{1}/{(2p)}<\g\le 1/2\,,
\end{array}\right.\\[4,5mm]
\hskip -5truemm
\label{WND}&W_{\textrm{N,D}}^{2\g,p}(\Omega)=
\left\{\!\!\begin{array}{l}
W^{2\g,p}(\Omega),\q\qq\qq\qq\qq\qq\qq\q\,
\textrm{if}\; 0<\g<{1}/{(2p)}\,,\\[1,5mm]
\big\{u\in W^{2\g,p}(\Omega):u=0 \;\;\textrm{on}\;
\partial\mbox{}B(0,R_1)
\big\},\q\;\;\,\textrm{if}\;\;{1}/{(2p)}<\g\le 1/2\,,
\end{array}\right.\\[4,5mm]
\hskip -5truemm
\label{WNN}& W_{\textrm{N,N}}^{2\g,p}(\Omega)=
 W^{2\g,p}(\Omega),\qq\qq\qq\qq\qq\qq\qq\q\,
\textrm{if}\;\;0<\g\le {1/2}\,.
\end{align}
\vskip 3truemm
\begin{remark}\label{embed}
\emph{ Observe that our choice $p\in (3,+\infty)$ implies the
embeddings (cf. \cite[Theorem 5.4]{AD})}
\begin{equation}\label{EMB}
 W^{1,p}(\Omega)\hookto C^{1-3/p}({\ov \Omega}),\qq
 W_{\!\emph{H,K}}^{2,p}(\Omega)\hookto C^{2-3/p}({\ov \Omega}).\qq
\end{equation}
\end{remark}
Assume also that $u_0$ satisfies the following conditions for some
positive constant $m$:\begin{eqnarray} \label{J0}
&J_0&\!\!\!\!\!(u_0)(r)\!:=\big|\Phi[\mathcal{C}u_0](r)\big|
\geqslant m\,,\qquad\;\forall\,r\in (R_1,R_2),\\[1,7mm]
\label{J1} &J_1&\!\!\!\!\!(u_0)\!:=\Psi[J(u_0)]\neq 0\,,
\end{eqnarray}
where we have set:
\begin{equation}\label{J}
J(u_0)(x)\!:=\!\bigg(\!\mathcal{B}u_0(x)-\frac{\Phi[\mathcal{B}u_0](|x|)}
{\Phi[\mathcal{C}u_0](|x|)}\mathcal{C}u_0(x)\!\bigg)\exp\!\bigg[
\int_{{\!|x|}}^{{R_2}}\!\frac{\Phi[\mathcal{B}u_0](\xi)}{\Phi[\mathcal{C}u_0]
(\xi)}d\xi\bigg]\,,\quad\forall\,x\in\Omega\,.
\end{equation}
\begin{remark}\emph{If $(\ref{ipotesiaij})$ and $(\ref{ipotesibijeci})$ hold,
then according to $(\ref{primasuPhiePsi})$ and
$(\ref{secondasuPhi})$ it follows that:
\begin{equation}
\Phi\big[J(u_0)\big](r)=\exp\!\bigg[\int_{{\!r}}^{{R_2}}\!\frac{\Phi[
\mathcal{B}u_0](\xi)}{\Phi[\mathcal{C}u_0](\xi)}d\xi\bigg]\Phi\bigg(\!\mathcal
{B}u_0- \frac{\Phi[\mathcal{B}u_0]}{\Phi[\mathcal{C}u_0]}
\mathcal{C}u_0\!\bigg)(r)=0\,,\quad\forall\,r\in (R_1,R_2)\,.
\end{equation}
\vskip -0.3truecm
\pn
This means that operator $\Psi$} cannot be chosen of the form
$\Psi\!=\! \Lambda\Phi$, where\, $\Lambda$ is a functional in $
L^p(R_1,R_2)^*$, \emph{otherwise condition $(\ref{J1})$ would be
not satisfied. In the explicit case, when $\Phi$ and $\Psi$ have
the integral representations $(\ref{Phi1})$ and $(\ref{Psi1})$,
this means that no function ${\psi}$ of the form
\,$\psi(x)={\psi}_1(|x|)\lambda\big(R_2\,{x}{|x|}^{-1}\big)$ is
allowed.}\end{remark}
\begin{remark}\emph{When operators $\Phi$ and $\Psi$ are defined by
$(\ref{Phi1})$, $(\ref{Psi1})$ conditions
$(\ref{J0})$, $(\ref{J1})$ can be rewritten as:}
\begin{align}
&\q\Big| \int_0^{\pi}\!\!\!\!\sin\!\theta
d\theta\!\!\int_0^{2\pi}\!\!\!\lambda(R_2x') \mathcal{C}u_0(rx')
d\varphi\,\Big|\!\geqslant m_1\,,
 \qq\forall\,r\in (R_1,R_2)\,,&\\[3mm]
&\q\bigg|\int_{\!{R_1}}^{{R_2}}\!\!\!\!r^2
dr\!\!\int_0^{\pi}\!\!\!\!\sin\!\theta d\theta\!\!
\int_0^{2\pi}\!\!\!\psi(rx')\bigg(\!\mathcal{B}u_0(rx')-\frac{\int_0^{\pi}\!
\!\sin\!\theta
d\theta\!\int_0^{2\pi}\!\!\lambda(R_2x')\mathcal{B}u_0(rx')
d\varphi}{\int_0^{\pi}\!\!\sin\!\theta
d\theta\!\int_0^{2\pi}\!\!\lambda(R_2x') \mathcal{C}u_0(rx')
d\varphi}\mathcal{C}u_0(rx')\!\bigg)&
\nonumber\\[2mm]
\label{bo}&\qq\qq\qq\;\;\;\times\exp\!\Bigg[\!\int_{{r}}^{{R_2}}
\!\!\!\frac{\;\int_0^{\pi}\!
\sin\!\theta
d\theta\!\int_0^{2\pi}\!\lambda(R_2x')\mathcal{B}u_0(\xi x')
d\varphi}{\int_0^{\pi}\!\sin\!\theta
d\theta\!\int_0^{2\pi}\!\lambda(R_2x') \mathcal{C}u_0(\xi x')
d\varphi}d\xi\Bigg] d\varphi \,\bigg|\geqslant m_2&
\end{align}
\emph{for some positive constants $m_1$ and $m_2$.}\end{remark}
Finally, we introduce the Banach spaces ${\mathcal{U}}^{\,s,p}(T)$,
${\mathcal{U}}_{\,\textrm{H,K}}^{\,s,p}(T)$, $\textrm{H,K}\in\{\textrm{D,N}\}$,
which are defined for any $s\in \mathbb{N}\backslash\{0\}$ by
\begin{equation}\label{Us}
\left\{\!\!\begin{array}{lll}
{\mathcal{U}}^{\,s,p}(T)\!\!\!&=&\!\!\!C^s\big([0,T];L^p(\Omega)\big)\cap
C^{s-1}\big([0,T];W^{2,p}(\Omega)\big)\,,\\[3mm]
{\mathcal{U}}_{\,\textrm{H,K}}^{\,s,p}(T)\!\!\!&=&\!\!\!C^s\big([0,T];L^p(\Omega)\big)\cap
C^{s-1}\big([0,T];W_{\!\textrm{H,K}}^{2,p}(\Omega)\big)\,.
\end{array}\right.
\end{equation}
Moreover, we list some further consistency conditions:
\begin{align}
\label{DD2}
&\qq(\textrm{C2,D,D})\q\qq\q\,
{v_0}(x)=0,\hskip 2,55truecm \forall\, x\in
\partial\mbox{}\Omega\,,&\\[2mm]
\label{DN2}
&\qq(\textrm{C2,D,N})\,\qquad\q\left\{\!\!\!
\begin{array}{lll} {v_0}(x)=0,
\hskip 2truecm &\forall&\!\!\!\! x\in \partial\mbox{}B(0,R_2)\;,
\\[2mm]
\displaystyle\frac{\partial v_0}{\partial\mbox{}\nu}(x)=0,
\hskip 2truecm
&\forall&\!\!\!\! x\in \partial\mbox{}B(0,R_1),\,\,
\end{array}\right.&\\[3mm]
\label{ND2}
&\qq(\textrm{C2,N,D})\quad\qquad\,\left\{\!\!\!
\begin{array}{lll}
\displaystyle\frac{\partial v_0}{\partial\mbox{}\nu}(x)=0,
\hskip 2truecm
&\forall&\!\!\!\! x\in \partial\mbox{}B(0,R_2),\,\\[3mm]
\, {v_0}(x)=0,\hskip 2truecm&\forall&\!\!\!\! x\in
\partial\mbox{}B(0,R_1),\,\,\,
\end{array}\right.&\\[3mm]
\label{NN2}
&\qq(\textrm{C2,N,N})\q\qq\;\;\,\;\frac{\partial
v_0}{\partial\mbox{}\nu}(x)=0,
\hskip 2,37truecm\forall\, x\in
\partial\mbox{}\Omega\,,&
\end{align}
\begin{eqnarray}
\label{PHIV1}&& \Phi[v_0](r)=D_tg_1(0,r)-\Phi[D_tu_1(0,\cdot)](r),\qq\forall r\in (R_1,R_2),\\[1,5mm]
\label{PSIV1}&& \Psi[v_0]=D_tg_2(0)-\Psi[D_tu_1(0,\cdot)]\,,
\end{eqnarray}
where \begin{equation}\label{v0}
v_0(x):\!={\mathcal{A}}u_0(x)+f(0,x)-D_tu_1(0,x)\,,\qquad
\forall\,x\in\Omega\,.
\end{equation}
\begin{theorem}\label{teoremaprincipe}
Let assumptions
$(\ref{unel}),\,(\ref{ipotesiaij})\!-\!(\ref{ipotesibijeci}),\,
(\ref{p3})\!-\!(\ref{primasuPsi})$ and condition $(\ref{HO})$ be
fulfilled. Moreover assume that the data enjoy the properties
$(\ref{richiestasuf})\!-\! (\ref{richiesteperg2})$ and  satisfy
inequalities $(\ref{J0}), (\ref{J1})$ and consistency conditions
$(\emph{C1,H,K})$ $($cf. $(\ref{DD1})-(\ref{NN1})$$)$, $(\emph{C2,H,K})$ as well
as $(\ref{1.18})$,
$(\ref{1.19})$, $(\ref{PHIV1})$, $(\ref{PSIV1})$.\\
 Then there exists
$T^{\ast}\in (0,T]$ such that the identification problem
$\emph\PIK$, $\emph{H,K}\in\{\emph{D,N}\}$,  admits a unique
solution $(u,k)\in
{\mathcal{U}}^{\,2,p}(T^{\ast})\times
C^{\beta}\big([0,T^{\ast}]; W^{1,p}(R_1,R_2)\big)$  depending
continuously on the data with respect to the norms related to the
Banach spaces in $(\ref{richiestasuf})\!-\!(\ref{richiesteperg2})$.\\
In the case of the specific operators $\Phi$, $\Psi$ defined by
$(\ref{Phi1}), \,(\ref{Psi1})$ the previous result is still true
if we assume that $\lambda\in C^1(\partial\mbox{}B(0,R_2))$ and
$\psi\in C^1(\overline{\Omega})$ with ${\psi}\!=\!0$ on the part
of $\partial\mbox{}\Om$ where the Dirichlet condition is
possibly prescribed.
\end{theorem}
\begin{lemma}\label{PHIPSI}
 When $\Phi$ and $\Psi$ are  defined by $(\ref{Phi1})$ and $(\ref{Psi1})$,
respectively, conditions
$(\ref{primasuPhiePsi})\!-\!(\ref{primasuPsi})$ are satisfied
under assumptions $(\ref{ipotesiaij}),\,(\ref{HO})$ on the
coefficients $a_{i,j}\;(i,j=1,2,3)$ and the hypotheses that
$\lambda \in C^1(\partial\mbox{}B(0,R_2))$  and $\psi\in
C^1(\overline{\Omega})$ with ${\psi}\!=\!0$ on the part of
$\partial\mbox{}\Om$ where the Dirichlet condition is possibly prescribed.
\end{lemma}
\begin{proof}
>From definitions $(\ref{Phi1}),\,(\ref{Psi1})$ it trivially
follows that conditions
$(\ref{primasuPhiePsi})\!-\!(\ref{terzasuPhi})$ are satisfied.
Hence we have only to prove that the decompositions
$(\ref{quartasuPhi}),\,
(\ref{primasuPsi})$ hold.\\
If the coefficients $a_{i,j}$, $i,j=1,2,3$, satisfy condition $(\ref{HO})$, then
the second-order differential operator $\mathcal{A}$ can be
represented, in spherical co-ordinates, by operator
$\widetilde{\mathcal{A}}$ defined by $(\ref{tildeA})$. Thus, taken
$w\in W_{\!\textrm{H,K}}^{2,p} (\Omega)$ with $p\in (3,+\infty)$, we can apply the
linear functional $\Phi$ defined in $(\ref{Phi1})$ to the
right-hand side of $(\ref{tildeA})$. From the well-known formulae
\begin{equation}
\label{Dr} \left\{\!\! \begin{array}{lll} D_{r}\!\!\! & =
&\!\!\!{\cos\!\varphi\sin\!\theta} D_{x_1}+\sin\!\varphi
\sin\!\theta D_{x_2}+\cos\!\theta D_{x_3}\,,\\[1,7mm]
D_{\varphi}\!\!\!  & = &\!\!\!{-r\sin\!\varphi\sin\!\theta}
D_{x_1}+
r\cos\!\varphi\sin\!\theta D_{x_2}\,, \\[1,7mm]
D_{\theta} \!\!\! & = &\!\!\!{r\cos\!\varphi\cos\!\theta}
D_{x_1}+r\sin\!\varphi\sin\!\theta D_{x_2}-r\sin\!\theta
D_{x_3}\,.
\end{array}\right.
\end{equation}
it follows that
$D_rw,\,{(D_{\varphi}w)}/(r\sin\!\theta),\,{(D_{\theta}w)}/{r}$
belong to $C^{1-3/p}([R_1,R_2]\times [0,2\pi]\times [0,\pi])$.
Hence, differentiating under the integral
sign, integrating by parts, using the periodicity with respect to
$\varphi$ of the functions $ g_j,k_j$, $j=1,2,3$, defined by
$(\ref{gj}),\,(\ref{kj,lj})$ and the membership of $\lambda$ in
$C^1(\partial\mbox{}B(0,R_2))$, we obtain:
\begin{eqnarray}\label{A_1w}
\Phi[\widetilde{\mathcal{A}}_1w](r)
 =   \int_0^{\pi}\!\!\!\!\sin\!\theta d\theta\!\!\int_0^{2\pi}\!\!\!
\lambda({R_2x}') D_r[k_1(r)D_r w(rx')] d\varphi=
\widetilde{\mathcal{A}}_1\Phi[w](r),\qquad\qquad\
\end{eqnarray}
$\widetilde{\mathcal{A}}_1$ being the differential operator
defined in $(\ref{tildeA1})$. Likewise we derive the formulae
\begin{align}
&\Phi \Big[ D_r
\Big(\frac{k_2(r,\varphi,\theta)
D_{\varphi}w}{r\sin\!\theta}\Big)\!\Big](r) =
-\!\int_0^{\pi}\!\!\!\! d\theta \!\!\int_0^{2\pi}\!\!\! D_r
w(rx')\frac{ D_{\varphi}[\lambda(R_2x')
k_2(r,\varphi,\theta)]}{r} d\varphi&\nonumber \\[3mm]
&\q+ \frac{1}{r}\int_0^{\pi}\!\!\!
\!d\theta\!\!\int_0^{2\pi}\! \!\! w(rx')\Big\{\frac{
D_{\varphi}[\lambda(R_2x')k_2(r,\varphi,\theta)]}{r}-
D_rD_{\varphi}[\lambda(R_2x')k_2(r,\varphi,\theta)] \Big\}
d\varphi,&
\end{align}
\begin{align}
&\Phi\Big[  D_r
\Big(\frac{k_3(r,\varphi,\theta) D_{\theta}w}{r}\Big)\!\Big](r) =
-\!\int_0^{\pi}\!\!\!\! d\theta\!\! \int_0^{2\pi}\!\!\!\! D_r
w(rx')\frac{ D_{\theta}[\lambda(R_2x')
k_3(r,\varphi,\theta)\!\sin\!\theta]}{r}\, d\varphi\;\;
\nonumber \\[3mm]
 &\q+ \frac{1}{r}\int_0^{\pi}\!\!\!\!d\theta\!\!\int_0^{2\pi}\!\!
\!\!\!  w(rx')\Big\{\!\frac{
D_{\theta}[\lambda(R_2x')k_3(r,\varphi,\theta)\!
\sin\!\theta]}{r}-D_rD_{\theta}[\lambda(R_2x')k_3(r,\varphi,\theta)\!\sin\!
\theta]\! \Big\} d\varphi,\\[4mm]
&\Phi \Big[ -\frac{ \sin\!\varphi}
{r\sin\!\theta} D_{\varphi}\Big({f_1(r,\varphi,\theta)D_rw}\!+\!
\frac{f_2(r,\varphi,\theta)D_{\varphi}w}{r\sin\!\theta}\!+\!
\frac{f_3(r,\varphi,\theta)D_{\theta}w}{r}\Big)\!\Big](r)&
\nonumber\\[3mm]
 &\q= \!\int_0^{\pi}\!\!\!\! d\theta\!\!\int_0^{2\pi}\Big[
{f_1(r,\varphi,\theta)D_rw(rx')}\!+\!\frac{f_2(r,\varphi,\theta)D_{\varphi}
w(rx')}{r\sin\!\theta}\!+\!\frac{f_3(r,\varphi,\theta)D_{\theta}w(rx')}{r}
\Big]& \nonumber \\[3mm]
 &\hskip 9,4truecm \times\frac{D_{\varphi}[{\lambda(R_2x')\!\sin\!\varphi}]}{r}\,
 d\varphi,&\\[4mm]
&\Phi \Big[ \frac{\cos\!\varphi}
{r\sin\!\theta } D_{\varphi}\Big({g_1(r,\varphi,\theta)D_rw}\!+\!
\frac{g_2(r,\varphi,\theta)D_{\varphi}w}{r\sin\!\theta}\!+\!
\frac{g_3(r,\varphi,\theta)D_{\theta}w}{r}\Big)\!\Big](r)\q\,\, \nonumber \\[3mm]
& \q =
-\!\int_0^{\pi}\!\!\!\!d\theta\!\!\int_0^{2\pi}\Big[
{g_1(r,\varphi,\theta)D_rw(rx')}\!+\!\frac{g_2(r,\varphi,\theta)D_{\varphi}
w(rx')}{r\sin\!\theta}\!+\!\frac{g_3(r,\varphi,\theta)
D_{\theta}w(rx')}{r}\Big]\!\! \nonumber \\[2,5mm]
 &  \hskip 9truecm \times \frac{D_{\varphi}[\lambda(R_2x')\!\cos\!\varphi]}{r}\,d\varphi,\\[4mm]
&\Phi  \Big[  -\frac{ \sin\!\theta}{r}
D_{\theta}\Big({h_1(r,\varphi,\theta)D_rw}\!+\!\frac{h_2(r,\varphi,\theta)
D_{\varphi}w}{r\sin\!\theta}\!+\!\frac{h_3(r,\varphi,\theta)D_{\theta}w}{r}
\Big)\!\Big](r)\nonumber \\[3mm]
&  \q =\!\int_0^{\pi}\!\!\!\!
d\theta\!\!\int_0^{2\pi}\Big[
h_1(r,\varphi,\theta)D_rw(rx')\!+\!\frac{h_2(r,\varphi,\theta)D_{\varphi}
w(rx')}{r\sin\!\theta}\!+\!\frac{h_3(r,\varphi,\theta)D_{\theta}w(rx')}{r}
\Big]\nonumber \\[2,5mm]
 &  \hskip 9truecm \times\frac{ D_{\theta}[\lambda(R_2x'){\sin}^2\theta]}{r}
 d\varphi,\\[3mm]
&\Phi  \Big[  \frac{ \cos\!\theta}{r}
D_{\theta}
\Big(l_1(r,\varphi,\theta)D_rw\!+\!\frac{l_2(r,\varphi,\theta)
D_{\varphi}w}{r\sin\!\theta}\!+\!\frac{l_3(r,\varphi,\theta)D_{\theta}w}{r}
\Big)\!\Big](r)\nonumber \\[3mm]
&
\q=-\!\int_0^{\pi}\!\!\!\!d\theta\!\!\int_0^{2\pi}\Big[
l_1(r,\varphi,\theta)D_rw(rx')\!+\!\frac{l_2(r,\varphi,\theta)D_{\varphi}
w(rx')}{r\sin\!\theta}\!+\!\frac{l_3(r,\varphi,\theta)D_{\theta}w(rx')}{r}
\Big]\quad \nonumber \\[3mm]
\label{17} &  \hskip 9truecm \times
\frac{D_{\theta}[\lambda(R_2x')\!\sin\!2\theta]}{2r}d\varphi.
\end{align}

After, rearranging the terms of $(\ref{A_1w})\!-\!(\ref{17})$ we
find that for every $w\in W_{\!\textrm{H,K}}^{2,p}(\Omega)$ with
$p\in (3,+\infty)$ the following equation holds: \vskip -0,1truecm
$$\Phi[\widetilde{\cal{A}}w](r)=\widetilde{\cal{A}}_1\Phi[w](r)+{\Phi}_1[w](r)\;,$$
\vskip -0,1truecm \pn where ${\Phi}_1$ is given by
\begin{align}
 &{\Phi} _1[ w](r)  =
\frac{1}{r}\displaystyle
\int_0^{\pi}\!\!\!\!d\theta\!\!\int_0^{2\pi}\!\!\!\!w(rx')\Big\{\displaystyle
\frac{ D_{\varphi}[\lambda(R_2x')k_2(r,\varphi,\theta)]}{r}+\frac{
D_{\theta}
[\lambda(R_2x')k_3(r,\varphi,\theta)\sin\!\theta]}{r}& \nonumber\\[2,5mm]
 &\q\; -  D_rD_{\varphi}[\lambda(R_2x')k_2(r,\varphi,\theta)]-
D_rD_{\theta}[\lambda(R_2x')k_3(r,\varphi,\theta)\!\sin\!\theta]\Big\}
d\varphi\,  +\int_0^{\pi}\!\!\!\!d\theta\!\!\int_0^{2\pi}\!\! D_r
w(rx')&
\nonumber\\[2,5mm]
 &\times\!\!\bigg\{\! \frac{f_1(r,\varphi,\theta)D_{\varphi}[
{\lambda(R_2x')\!\sin\!\varphi}]}{r}-\frac{g_1(r,\varphi,\theta)
D_{\varphi}[{\lambda(R_2x')\!\cos\!\varphi}]}{r}
+\frac{h_1(r,\varphi,
\theta)D_{\theta}[{\lambda(R_2x'){\sin}^2\theta}]}{r}&
\nonumber\\[2,5mm]
 & \,-\frac{l_1(r,\varphi,\theta)D_{\theta}[\lambda(R_2x'){\sin\!2
\theta}] }{2r}- \frac{
D_{\varphi}[\lambda(R_2x')k_2(r,\varphi,\theta)]} {r}-\frac{
D_{\theta}[\lambda(R_2x')k_3(r,\varphi,\theta)\!\sin\!\theta]}
{r}
\bigg\}d\varphi\,& \nonumber\\[2,5mm]
& \q\;+\! \int_0^{\pi}\!\!\!\!  d\theta\!\!\int_0^{2\pi}
\frac{D_{\varphi}w(rx')}{r\sin\!\theta}\bigg\{\frac{f_2(r,\varphi,\theta)
D_{\varphi}[{\lambda(R_2x')\!\sin\!\varphi}]}{r}-\frac{g_2(r,\varphi,\theta)
D_{\varphi}[{\lambda(R_2x')\!\cos\!\varphi}]}{r}&\nonumber\\[2,5mm]
 &\qq\qq +\frac{h_2(r,\varphi,\theta)D_{\theta}[{\lambda(R_2x')
{\sin}^2\theta}]}{r}-\frac{l_2(r,\varphi,\theta)D_{\theta}[\lambda(R_2x')\!
\sin\!2\theta] }{2r}\bigg\}d\varphi&\nonumber\\[2,5mm]
 &\q\;+\!\int_0^{\pi}\!\!\!\!d\theta\!\!\int_0^{2\pi}\frac{D_{\theta}
w(rx')}{r}\bigg\{\frac{f_3(r,\varphi,\theta)D_{\varphi}[{\lambda(R_2x')\!
\sin\!\varphi}]}{r}-\frac{g_3(r,\varphi,\theta)D_{\varphi}[{\lambda(R_2x')\!
\cos\!\varphi}]}{r}& \nonumber\\[2,5mm]
\label{phi1} & \qq\qq+\,\frac{h_3(r,\varphi,\theta)
D_{\theta}[{\lambda(R_2x'){\sin}^2\theta}]}{r}-\frac{l_3(r,\varphi,\theta)
D_{\theta}[\lambda(R_2x')\!\sin\!2\theta] }{2r}\bigg\} d\varphi\,.&
\end{align}
We now prove that $\Phi_1$ belongs to
$\mathcal{L}\big(W^{1,p}({\Om}); L^p(R_1,R_2)\big)$, i.e. it
is bounded from $W^{1,p}({\Om})$ to
$L^p(R_1,R_2)$.\\
By virtue of assumption $(\ref{ipotesiaij})$ on the coefficients $
 a_{i,j}$, $i,j=1,2,3,$ from $(\ref{fj})\!-\! (\ref{hj})$, $(\ref{kj,lj})$ we
deduce that $f_j,g_j,h_j,k_j,l_j$ belong to $
W^{2,\infty}((R_1,R_2)\times (0,2\pi)\times (0,\pi))$, $j=1,2,3$.
Then, since  $\lambda\in C^1(\partial\mbox{}B(0,R_2))$ and
$0\!<\!R_{2}^{-1}\!<\!r^{-1}\!<\!R_{1}^{-1}$, using formulae
$(\ref{Dr})$ we can  prove that the following functions
 $$\begin{array}{l}
\qquad\qquad\quad
D_rD_{\theta}[\lambda(R_2x')k_3(r,\varphi,\theta)\!
\sin\!\theta]\;,\; D_rD_{\varphi}[\lambda(R_2x')k_2(r,\varphi,\theta)]\;,\\ \\
\quad\displaystyle\frac{
D_{\varphi}[\lambda(R_2x')k_2(r,\varphi,\theta)]}{r} \:,\;
\displaystyle\frac{
D_{\theta}[\lambda(R_2x')k_3(r,\varphi,\theta)\!
\sin\!\theta]}{r}\;,\;\displaystyle
\frac{f_j(r,\varphi,\theta)D_{\varphi}
[{\lambda(R_2x')\!\sin\!\varphi}]}{r}\;,
\end{array}$$
\begin{equation}
\displaystyle\frac{g_j(r,\varphi,\theta)D_{\varphi}[{\lambda(R_2x')\!
\cos\!\varphi}]}{r}\;,\; \displaystyle\frac{h_j(r,\varphi,\theta)
D_{\theta}[{\lambda(R_2x')]{\sin}^2\theta}]}{r}\;,\; \displaystyle
\frac{l_j(r,\varphi,\theta)D_{\theta}[\lambda(R_2x')\!\sin\!2\theta]
}{2r} \nonumber
\end{equation}
belong to $ L^{\infty}((R_1,R_2)\times (0,2\pi)\times (0,\pi))$
and their $L^{\infty}$-norms are bounded from above by
$C\|\lambda\|_{C^1(\partial\mbox{}B(0,R_2))}$, $C$ being a positive
constant depending on $R_1,R_2,\max_{\,_{_{i,j=1,2,3}}}
\|a_{i,j}\!\|_{W^{2,\infty}({\Omega})}$,  only. \pn Observe now
that for any pair of functions $f\in C(\overline{\Omega})$ and
$v\in L^{p}(\Omega)$ we have
\begin{eqnarray}
\int_0^{\pi}\!\!\!\!d\theta\!\!\int_0^{2\pi}\!\!|v(rx')f(rx')|d\varphi\!\!\!
&\leqslant&\!\!\!
{\|f\|}_{C(\ov\Om)}\!{\bigg[\!\int_0^{\pi}\!\!\!\!\sin\!\theta
d\theta\!\! \int_0^{2\pi}\!\!\!|v(rx'){|}^{p}d\varphi\bigg]}^{\!
{1}/{p}}{\bigg[2\pi\!\!\int_0^{\pi}\!\!\!{(\sin\!\theta)}^{\!^{\!-\frac{1}{p-1}}}
d\theta\bigg]}^{\!{p}/{(p-1)}}\nonumber\\
\label{perche'p>3}
\end{eqnarray}
Since the right-hand side in $(\ref{perche'p>3})$ is in
$L^{p}(R_1,R_2)$ when $p\in (3,+\infty)$, applying H\"older's inequality to the
right-hand side of $(\ref{phi1})$ we find
$\|\Phi_1[w]\|_{L^p(R_1,R_2)} \leqslant
C\|w\|_{W^{1,p}({\Omega})}$. Consequently decomposition
$(\ref{quartasuPhi})$ holds.

Let now  $\Psi\in {L^p(\Omega)}^{\ast}$ be the functional defined
in $(\ref{Psi1})$. Analogously to what we have done for $\Phi$, we
apply $\Psi$ to both sides in  $(\ref{tildeA})$. Performing
computations similar to those made above and using the assumption
that $\psi_{|_{\G}}\!=\!0$, when the Dirichlet condition
is prescribed on $\G\subset\partial\mbox{}\Om$, we obtain the
equation:
$$\Psi[\widetilde{\cal{A}}w]= {\Psi}_1[w]\,,\qq\forall\,w\in W_{\textrm{H,K}}^{2,p}(\Om)\,,$$
where
\begin{align}
&{\Psi}_1[ w]=
 -\!\int_{R_1}^{R_2}\!\!\!\!rdr\!\!\int_0^{\pi}\!\!\!\!d\theta\!\!
\int_0^{2\pi}
\!\!\!D_rw(rx')\bigg\{\frac{k_1(r)D_r[r^2\psi(rx')]\!\sin\!\theta}{r}
-f_1(r,\varphi,\theta)
D_{\varphi}[\psi(rx')\!\sin\!\varphi]&\nonumber\\[3mm]
&+ g_1 (r,\varphi,\theta)D_{\varphi}[\psi(rx')
\!\cos\!\varphi]-h_1(r,\varphi,\theta)D_{\theta}[{\psi(rx')
{\sin}^2\theta}]+\frac{
l_1(r,\varphi,\theta)D_{\theta}[\psi(rx')\!
\sin\!2\theta]}{2}\bigg\}d\varphi&\nonumber\\[3mm]
&\;-\int_{R_1}^{R_2}\!\!\!\!rdr\!\!\int_0^{\pi}\!\!\!\!
d\theta\!\!\int_0^{2\pi}\!\frac{D_{\varphi}w(rx')}{r\sin\!\theta}
\bigg\{\!\frac{k_2(r,\varphi,\theta)D_r[r^2\psi(rx')]\!\sin\!\theta}{r}
-f_2(r,\varphi,\theta)
D_{\varphi}[\psi(rx')\!\sin\!\varphi]& \nonumber \\[3mm]
&+g_2(r,\varphi,\theta)D_{\varphi}[\psi(rx')\!
\cos\!\varphi]-h_2(r,\varphi,\theta)D_{\theta}[{\psi(rx')
{\sin}^2\theta}]+\frac{
l_2(r,\varphi,\theta)D_{\theta}[\psi(rx')
\sin\!2\theta]}{2}  \bigg\}d\varphi&\nonumber\\[3mm]
&\; -
\int_{R_1}^{R_2}\!\!\!\!rdr\!\!\int_0^{\pi}\!\!\!\!
d\theta\!\!\int_0^{2\pi}\frac{D_{\theta}w(rx')}{r}\bigg\{\!
\frac{k_3(r,\varphi, \theta) D_r[r^2\psi(rx')]\!\sin\!\theta}{r}
-f_3(r,\varphi,\theta)D_{\varphi}[
\psi(rx')\!\sin\!\varphi]&\nonumber\\[3mm]
&+ g_3 (r,\varphi,\theta)D_{\varphi}[\psi(rx')\!
\cos\!\varphi]-h_3(r,\varphi,\theta)D_{\theta}[{\psi(rx')
{\sin}^2\theta}]+\frac{
l_3(r,\varphi,\theta)D_{\theta}[\psi(rx')\!
\sin\!2\theta]}{2}\bigg\}d\varphi\,.&\nonumber\\
\label{ps1}& &
\end{align}
Since the functions $f_j,g_j,h_j,k_j,l_j$, $j\!=\!1,2,3$, belong to $
W^{2,\infty}({\Omega})$ and $\psi\!\in\! C^1(\overline{\Omega})$,
using an estimate similar to $(\ref{perche'p>3})$, it easily
follows that ${\Psi}_1\in {W_{\textrm{H,K}}^{1,p} ({\Omega})}^{\ast}$. Hence
decomposition $(\ref{primasuPsi})$ also holds. This completes
the proof.
\end{proof}
\begin{remark}\label{nosphere}
\emph{The reason why we have restricted ourselves to investigating
the identification problem $\PIK$ in the spherical corona
$\Omega\!=\!\{x\in\mathbb{R}^3\!:\!R_1<|x| <R_2\}$, $0<R_1<R_2$,
instead of the simpler ball ${\Omega}_1\!
=\!\{x\in\mathbb{R}^3\!:\!|x|<R\}$, $R>0$, is due to the
representation $(\ref{phi1})$ of the functional $\Phi_1$. Indeed,
the function
appearing in the right-hand side of $(\ref{phi1})$ might not
belong to $L^p(0,\ve)$ for any $\ve \in (0,R_2)$ when dealing with
general coefficients $a_{i,j}\in W^{2,\infty}(\Omega_1)$.
 This would imply
$\Phi_1\not \in \mathcal{L}\big(W^{1,p}({\Omega}_1);L^p(0,R)\big)$ and
would prevent us from applying known abstract results.}
\end{remark}
\section{An equivalence result in the concrete case}
\setcounter{equation}{0} In this section we prove an equivalence
theorem which will be the starting point to reduce our problem to
the same abstract integral fixed-point system
studied in \cite{CL}.\\
Let us suppose that
$(u,k)\in{\mathcal{U}}^{\,2,p}(T)\times
C^{\beta}\big([0,T];W^{1,p}(R_1,R_2)\big)$ is a solution to the
identification
 problem $\PIK$. Let us now introduce the following new unknown function
\begin{eqnarray}\label{v}
&& v(t,x)=D_tu(t,x)-D_t{u}_1(t,x) \q \Longleftrightarrow \qq \nonumber
\\[2mm]
&& u(t,x)=u_1(t,x)-u_1(0,x)+u_0(x) +\int_0^t\!v(s,x)ds.
\end{eqnarray}
Then from $(\ref{problem})$ it follows that the pair
$(v,k)\in {\mathcal{U}}_{\,\textrm{H,K}}^{\,1,p}(T)\times
C^{\beta}\big([0,T];W^{1,p}(R_1,R_2)\big)$ solves the
identification problem
\begin{eqnarray}
\label{problem1}
D_tv\!\!\!\!\!&(t,x)&\!\!\!\!=\mathcal{A}v(t,x)+\int_0^t\!
k(t-s,|x|)
\big[\mathcal{B}v(s,x)+\mathcal{B}D_{t}u_1(s,x)\big]ds+k(t,|x|)
\mathcal{B}u_0(x)\nonumber\\[1,2mm]
&&\;+\!\int_0^t\! D_{|x|}k(t-s,|x|)\big[\mathcal{C}v(s,x)+
\mathcal{C}D_t{u}_1(s,x)\big]ds+D_{|x|}k(t,|x|)\mathcal{C}u_0(x)
\nonumber\\[1,4mm]
& &\;+\,\mathcal{A}D_t{u}_1(t,x)-D_t^{2}{u}_1(t,x)
+D_tf(t,x),\qq\;\forall\,(t,x)\in[0,T]\times\Omega\,,
\end{eqnarray}
\begin{align}
\label{3.3}&v(0,x)={\mathcal{A}}u_0(x)+f(0,x)-D_tu_1(0,x)=:v_0(x),\qquad
\forall\,x\in\Omega\,,&\\[1,5mm]
&v\; \textrm{satisfies the homogeneous boundary
conditions (H,K)}\,,&\\[2mm]
\label{PHIV}&\Phi[v(t,\cdot)](r)=D_tg_1(t,r)-\Phi[D_tu_1(t,\cdot)](r)\,,\q\;\forall\,(t,r)\in[0,T]\times(R_1,R_2)\,,&
\\[1,7mm]
\label{PSIV}&\Psi[v(t,\cdot)]=D_tg_2(t)-\Psi[D_tu_1(t,\cdot)]\,,\qq\q
\forall\, t\in[0,T].&
\end{align}
The consistency conditions related to problem $(\ref{problem1})\!-\!
(\ref{PSIV})$ can be deduced as in section 1
with $(u_0,u_1, g_1, g_2)$ replaced by $(v_0, 0, D_tg_1-\Phi[D_tu_1],
D_tg_2-\Psi[D_tu_1])$ and they are explicitly given by $(\ref{DD2})\!-\!
(\ref{PSIV1})$.\\
 Using assumptions
$(\ref{primasuPhiePsi})-(\ref{primasuPsi})$ and applying the
functionals $\Phi$, $\Psi$ to both sides of $(\ref{problem1})$, it
easy to check that the radial kernel $k$ satisfies the two
following equations:
\begin{align}
\label{primasuk}
&D_rk(t,r)\Phi[\mathcal{C}u_0](r)+k(t,r)
\Phi[\mathcal{B}u_0](r)=N_1^{0}(u_1,g_1,f)(t,r)
+\Phi[\widetilde{N}_{1}(v,k)(t,\cdot)](r)&
\nonumber\\[1,6mm]
&\hskip 6truecm -\,{\Phi}_1[v(t,\cdot)](r),\!\qq
 \forall\,(t,r)\in [0,T]\!\times\!(R_1,R_2),&\\[2,8mm]
&\Psi[D_rk(t,\cdot)\mathcal{C}u_0+k(t,\cdot)\mathcal{B}u_0]=
N_2^{0}(u_1,g_2,f)(t)+
\Psi[\widetilde{N}_{1}(v,k)(t,\cdot)]
-{\Psi}_1[v(t,\cdot)],&\nonumber\\[1mm]
&\hskip 11,5truecm\label{secondasuk}\,\;\forall\,t\in[0,T],&
\end{align}
where the operators $\widetilde{N}_1,\,N_1^0,\,N_2^0$ are defined,
respectively, by
\begin{align}
&\widetilde{N}_1(v,k)(t,r)=-\int_0^t\!
k(t-s,|x|)
\big[\mathcal{B}v(s,x)+\mathcal{B}D_{t}u_1(s,x)\big]ds&
\nonumber\\[1,5mm]
\label{N1tilde}&\hskip 1,7truecm -\!\int_0^t\!\!
D_{|x|}k(t-s,|x|)
\big[\mathcal{C}v(s,x)+\mathcal{C}D_t{u}_1(s,x)\big]ds\,,
\q\forall\,(t,x)
\in[0,T]\!\times\!\Omega,&
\end{align}
\begin{eqnarray}
\hskip -7truemm
 N_1^{0}(u_1,g_1,f)(t,r)\!\!\!&=&\!\!\! D_t^2g_1(t,r)
-D_t{\widetilde{\mathcal{A}}}_1g_1(t,r)-{\Phi}_1[D_tu_1(t,\cdot)](r)-\Phi[D_tf(t,\cdot)](r)\,,
\nonumber\\[2mm]
\hskip -7truemm&&\label{N10}\hskip 5,3truecm
\forall\,(t,r)\in[0,T]\!
\times\!(R_1,R_2),\\[3mm]
\hskip -7truemm\label{N20}N_2^{0}(u_1,g_2,f)(t)\!\!\!&
=&\!\!\! D_t^2g_2(t)-{\Psi}_1[D_tu_1(t,\cdot)]
-{\Psi}[D_t f(t,\cdot)]\,,\qq\,\forall\,t\in[0,T].
\end{eqnarray}
\begin{remark}\emph{It can be easily checked that if $(v,k)
\in{\mathcal{U}}_{\,\textrm{H,K}}^{\,1,p}(T)\times
C^{\beta}\big([0,T];W^{1,p}(R_1,R_2)\big)$ solves the
identification problem
 $(\ref{problem1})\!-\!(\ref{secondasuk})$
then, taking advantage of the consistency conditions
$(\ref{DD2})\!-\!(\ref{PSIV1})$, the function
$u\in{\mathcal{U}}^{\,2,p}(T)$ defined in
$(\ref{v})$ is a solution to the problem  $\PIK$.}
\end{remark}
>From $(\ref{primasuk}),\,(\ref{secondasuk})$ it turns out that the
initial value $k(0,\cdot)$ must satisfy the following equations:
\begin{align}
\label{primasuk0}
&D_rk(0,r)\Phi[\mathcal{C}u_0](r)\!+\!
k(0,r)\Phi[\mathcal{B}u_0](r)=
N_1^{0}(u_1,g_1,f)(0,r)-{\Phi}_1[v_0](r),
\q \forall\,r\in (R_1,R_2),\nonumber\\[1mm]\\[2mm]
\label{secondasuk0}
&\Psi[D_rk(0,\cdot)\mathcal{C}u_0+k(0,
\cdot)\mathcal{B}u_0]=N_2^{0}(u_1,g_2,f)(0)-
{\Psi}_1[v_0].
\end{align}
Let
\begin{equation}\label{l1}
\wtil{l}_1(r)\!:=N_1^{0}(u_1,g_1,f)(0,r)-\!{\Phi}_1[v_0](r),\qquad\forall\;r\in
(R_1,R_2).
\end{equation}
Then using condition $(\ref{J0})$ and integrating the first-order
differential equation $(\ref{primasuk0})$ we obtain the following
general integral depending on an arbitrary constant $C$:
\begin{equation}\label{intk0}
k(0,r)=C\exp\!\bigg[\!\int_{{\!r}}^{{R_2}}\frac{\Phi[\mathcal{B}u_0](\xi)}
{\Phi[\mathcal{C}u_0](\xi)}d\xi\bigg]+\int_{\!R_2}^{r}\!\!\exp\!\bigg[\!
\int_{\!r}^{\eta}\frac{\Phi[\mathcal{B}u_0](\xi)}{\Phi[\mathcal{C}u_0](\xi)}
d\xi\bigg]\frac{\wtil{l}_1(\eta)}{\Phi[\mathcal{C}u_0](\eta)}d\eta\,.
\end{equation}
Substituting this representation of $k(0,\cdot)$ into
$(\ref{secondasuk0})$, we can compute the constant $C$:
\begin{eqnarray}
\label{Ck0}C\!\!\!\!&=&\!\!\!\!{[J_1(u_0)]}^{-1}\Big\{\Psi[\wtil{l}_2]
+N_2^{0}(u_1,g_2,f)(0)-{\Psi}_1[v_0]\Big\},\qquad
\end{eqnarray}
where $J_1(u_0)$ and $\wtil{l}_2$ are defined, respectively, by
$(\ref{J1})$ and the following formula:
\begin{eqnarray}
\wtil{l}_2(x)\!\!\!\!&:= &\!\!\!\!\mathcal{C}u_0(x)\bigg\{
\frac{\wtil{l}_1(|x|)}
{\Phi[\mathcal{C}u_0](|x|)}-\frac{\Phi[\mathcal{B}u_0](|x|)}
{\Phi[\mathcal{C}u_0](|x|)}\int_{\!R_2}^{|x|}\!\!\!\exp\!
\bigg[\int_{\!|x|}^{\eta}\frac{\Phi[\mathcal{B}u_0](\xi)}
{\Phi[\mathcal{C}u_0]
(\xi)}d\xi\bigg]\frac{\wtil{l}_1(\eta)}{\Phi[\mathcal{C}u_0](\eta)}
d\eta\bigg\}\nonumber\\\nonumber\\
\label{l2}&
&\!\!\!\!+\mathcal{B}u_0(x)\int_{\!R_2}^{|x|}\!\!\!\exp\!\bigg[
\int_{\!|x|}^{\eta}\frac{\Phi[\mathcal{B}u_0](\xi)}
{\Phi[\mathcal{C}u_0](\xi)}d\xi\bigg]\frac{\wtil{l}_1(\eta)}
{\Phi[\mathcal{C}u_0](\eta)}d\eta,\quad\qquad\forall
\;x\in\Omega\,.
\end{eqnarray}
Then, substituting $(\ref{Ck0})$ into $(\ref{intk0})$, we find
that the initial value $k(0,\cdot)$ is given by
\begin{eqnarray}
\label{k01} k(0,r)\!\!\!\!&=&\!\!\!\![J_1(u_0){]}^{-1}\Big\{
\Psi[\wtil{l}_2]
+N_2^{0}(u_1,g_2,f)(0)-{\Psi}_1[v_0]\!\Big\}
\exp\!\bigg[\int_{\!r}^{R_2}\frac{\Phi[
\mathcal{B}u_0](\xi)}{\Phi[\mathcal{C}u_0](\xi)}d\xi\bigg]
\nonumber\\[1,2mm]
&&\!\!\!\!+\int_{\!R_2}^{{r}}
\!\!\!\exp\!\bigg[\int_{\!r}^{\eta}\!\frac{\Phi[\mathcal{B}u_0](\xi)}
{\Phi[\mathcal{C}u_0](\xi)}d\xi\bigg]\frac{\wtil{l}_1(\eta)}
{\Phi[\mathcal{C}u_0](\eta)}d\eta:= k_0(r),\q\; \forall\;r\in (R_1,R_2)\,.
\end{eqnarray}
Now we introduce the two new unknown functions
\begin{equation}\label{q}
\hskip 0,5truecm h(t)=k(t,R_2),\qquad q(t,r)=D_rk(t,r),\;\quad\forall\;
(t,r)\in[0,T]\times(R_1,R_2).
\end{equation}
and express $k$ in terms of $h$ and $q$:
\begin{equation}\label{kkk}
k(t,r)=h(t)-\int_{r}^{R_2}\!\!\!\!q(t,\xi)d\xi:=h(t)-Eq(t,r)\,,
\;\quad\forall\;
(t,r)\in[0,T]\times(R_1,R_2).
\end{equation}
Of course, $(\ref{q})$ and $(\ref{kkk})$ imply the initial conditions
\begin{equation}
h(0)=k_0(R_2),\qq q(0,r)=k_0'(r),\qquad \forall\;r\in [R_1,R_2].
\end{equation}
Using $(\ref{kkk})$, we solve $(\ref{primasuk})$,
$(\ref{secondasuk})$ for the pair $(h,q)$. From definition
$(\ref{N1tilde})$ we deduce the following representation for
operator $\widetilde{N}_1$:
\begin{eqnarray}
\wtil{N}_1(v,h-Eq)(t,|x|)\!\!\!&:=&\!\!\!\!\!-\!\int_0^t\!\big[
h(t-s)-Eq(t-s,|x|)
\big]\big[\mathcal{B}v(s,x)+\mathcal{B}D_{t}u_1(s,x)\big]ds\nonumber\\
&&\!\!\!\!\!-\!\int_0^t\!\!q(t-s,|x|)\big[\mathcal{C}v(s,x)+
\mathcal{C}
D_t{u}_1(s,x)\big]ds\nonumber\\[2mm]
\label{N1}&:=&\!\!\!N_1(v,h,q)(t,|x|),\qquad
\forall\,(t,x)\in[0,T]\!\times\!\Omega\,.
\end{eqnarray}
Moreover system $(\ref{primasuk})$,\,$(\ref{secondasuk})$ changes
into
\begin{eqnarray}
q(t,r)\Phi[\mathcal{C}u_0](r)\!\!\!\!&-&\!\!\!\!Eq(t,r)\Phi[
\mathcal{B}u_0](r)=N_1^{0}(u_1,g_1,f)(t,r)-h(t)
\Phi[\mathcal{B}u_0](r)-{\Phi}_1[v(t,\cdot)](r)
\nonumber\\[1,6mm]
\label{qEq} & & \qq\q\, +\, \Phi[N_1(v,h,q)(t,\cdot)](r),
\quad\;\forall\,(t,r)\in[0,T]\!\times\! (R_1,R_2),\\[2,8mm]
\Psi\!\big[q(t,\cdot)\mathcal{C}u_0+[\!\!\!\!\!&h&\!\!\!\!\!\!(t)
-Eq(t,\cdot)]
\mathcal{B}u_0\big]=N_2^{0}(u_1,g_2,f)(t)+\,\Psi[{N_{1}}(v,h,q)(t,\cdot)]
-{\Psi}_1[v(t,\cdot)],
\nonumber\\[1,6mm]
&&\hskip 8truecm \qquad \forall\,t\in[0,T]. \label{hEq}
\end{eqnarray}
First we consider the integral equation
\begin{equation}\label{g}
q(t,r)\Phi[\mathcal{C}u_0](r)-Eq(t,r)\Phi[\mathcal{B}u_0](r)=g(t,x),\;\;\;
\quad\forall\;(t,r)\in[0,T]\!\times\! (R_1,R_2),
\end{equation}
where $g\in L^1\big((0,T)\!\times\! (R_1,R_2)\big)$ is an arbitrary
given function.\\
Since $u_0$ satisfies condition $(\ref{J0})$ and $(\ref{kkk})$
implies $Eq(t,R_2)=0$ for all $t\in [0,T]$  and
$D_rEq(t,r)\!=\!-q(t,r)$ for all $(t,r)\in[0,T]\times (R_1,R_2)$,
the solution to the
differential equation $(\ref{g})$ is given by
\begin{equation}\label{EL}
Eq(t,r)=Lg(t,r)\,,
\end{equation}
where operator $L$ is defined by the formula
\begin{equation}\label{L}
Lg(t,r)\!:=\int_{r}^{R_2}\!\exp\!\bigg[\int_{\!r}^{\eta}\frac{\Phi[
\mathcal{B}u_0](\xi)}{\Phi[\mathcal{C}u_0](\xi)}d\xi\bigg]\frac{g(t,\eta)}
{\Phi[\mathcal{C}u_0](\eta)}d\eta\,.
\end{equation}
Hence, using $(\ref{EL})$ and the relation
$D_rEq(t,r)\!=\!-q(t,r)$, from $(\ref{L})$ we obtain the following
representation formula for $q$:
\begin{equation}\label{q1}
\qquad
q(t,r)=\frac{1}{\Phi[\mathcal{C}u_0](r)}\big[I+\Phi[\mathcal{B}u_0](r)L
\big]g(t,r),\qquad\forall\;(t,r)\in[0,T]\!\times\! (R_1,R_2),
\end{equation}
where $I$ denotes the identity operator.\\
>From $(\ref{qEq})$ we get:
$$g(t,r)=N_1^{0}(u_1,g_1,f)(t,r)-h(t)\Phi[\mathcal{B}u_0](r)+\Phi[{N_{1}}
(v,h,q)(t,\cdot)](r)-{\Phi}_1[v(t,\cdot)](r).$$
Therefore
substituting into $(\ref{q1})$ we find:
\begin{eqnarray}
q(t,r)\!\!\!&=&\!\!\!-h(t)\frac{\Phi[\mathcal{B}u_0](r)}
{\Phi[\mathcal{C}u_0]
(r)}\big[1+L\Phi[\mathcal{B}u_0](r)\big]+N_3^0(u_0,u_1,g_1,f)(t,r)
+N_2(v,h,q)
(t,r),\nonumber\\
\label{q2}&&\qquad\qquad\qquad\qquad\qquad\qquad\qquad\qquad\qquad\forall\;
(t,r)\in[0,T]\!\times\! (R_1,R_2),
\end{eqnarray}
where we have set:
\begin{eqnarray}
\q\;
N_2(v,h,q)(t,r)\!\!\!&=&\!\!\!\frac{1}{\Phi[\mathcal{C}u_0](r)}\big[I+\Phi[
\mathcal{B}u_0](r)L\big]\big\{\Phi[{N_{1}}(v,h,q)(t,\cdot)](r)-{\Phi}_1
[v(t,\cdot)](r)\big\}\quad\nonumber\\[2mm]
\label{N2}\!\!\!& =:\!\!\!& J_3(u_0)(r)
\big\{\Phi[{N_{1}}(v,h,q)(t,\cdot)](r)-{\Phi}_1[v(t,\cdot)](r)\big\},
\end{eqnarray}
\begin{eqnarray}
N_3^0(u_0,u_1,g_1,f)(t,r)\!\!\!&=\!\!\!&\frac{1}{\Phi[\mathcal{C}u_0](r)}
\big[I+\Phi[\mathcal{B}u_0](r)L\big]N_1^{0}(u_1,g_1,f)(t,r)\qquad\qquad
\qquad\nonumber\\[2mm]
\label{N30}\!\!\!&=:\!\!\!&\!J_3(u_0)(r)N_1^{0}(u_1,g_1,f)(t,r)\,.
\qq\q\,
\end{eqnarray}
Observing that $(\ref{L})$ implies
\begin{equation}\label{I+L}
1+L\Phi[\mathcal{B}u_0](r)=\exp\!\bigg[\!\int_{\!{r}}^{\!{R_2}}
\frac{\Phi[\mathcal{B}u_0](\xi)}{\Phi[\mathcal{C}u_0](\xi)}d\xi\bigg],\qquad
\forall\,r\in (R_1,R_2),
\end{equation}
and substituting this expression into $(\ref{q2})$, from
$(\ref{hEq})$ it is easy to check that $h$ solves the following
equation:
\begin{eqnarray}h(t)J_1(u_0)\!\!\!&=&\!\!\!\!N_0(u_0,u_1,g_1,g_2,f)(t)
+\Psi[{N_{1}}(v,h,q)(t,\cdot)]\!-\!\Psi[N_2(v,h,q)(t,\cdot)\mathcal{C}u_0]
\nonumber\\[2mm]
\label{hh}&&\!\!\!\!+ \Psi\big[E\big(N_2(v,h,q)(t,\cdot)\big)
\mathcal{B}u_0\big]\!-\!{\Psi}_1[v(t,\cdot)]\,,\!\quad\quad\quad\forall\,
t\in[0,T]\,,
\end{eqnarray}
where $J_1(u_0)$ and $N_0(u_0,u_1,g_1,g_2,f)$ are defined,
respectively, by $(\ref{J1})$ and
\begin{eqnarray}
N_0(u_0,u_1,g_1,g_2,f)(t)\!:=N_2^{0}(u_1,g_2,f)(t)-
\Psi[N_3^0(u_0,u_1,g_1,f)(t,\cdot)\mathcal{C}u_0]\;\qquad\nonumber\\[1,8mm]
\label{N0}\qq-\Psi\big[E\big(N_3^0(u_0,u_1,g_1,f)(t,\cdot)\big)
\mathcal{B}u_0\big],\qquad \forall\,t\in[0,T].
\end{eqnarray}
Hence, from $(\ref{hh})$ and (\ref{J1}) we conclude that $h$ solves
the following fixed-point equation:
\begin{equation}\label{hhh}
h(t)= h_0(t)+N_3(v,h,q)(t),\qq \forall\;t\in[0,T],
\end{equation}
where we have set:
\begin{eqnarray}
\label{h0}
&& h_0(t)\!:={\big[J_1(u_0)\big]}^{-1}N_0(u_0,u_1,g_1,g_2,f)(t),\\[1,3mm]
\label{N3} && N_3(v,h,q)(t)\!:={\big[J_1(u_0)\big]}^{-1}
\Big\{\Psi[{N_{1}}(v,h,q)(t,\cdot)]\!-\!\Psi[N_2(v,h,q)(t,\cdot)
\mathcal{C}u_0]\nonumber\\
& &\hskip 3truecm +\Psi\big[E\big(N_2(v,h,q)(t,\cdot)\big)
\mathcal{B}u_0\big]\!-\!{\Psi}_1[v(t,\cdot)]\Big\}.
\end{eqnarray}
So, using again $(\ref{I+L})$ and replacing the right-hand side of
$(\ref{hhh})$ into $(\ref{q2})$, we conclude that  $q$ satisfies
the following fixed-point equation
\begin{eqnarray}
 q(t,r)=q_0(t,r)+J_2(u_0)(r)N_3(v,h,q)(t)
+N_2(v,h,q)(t,r),\nonumber\\[1,5mm]
\label{q3}\hskip 6truecm\forall\,(t,r)\in [0,T]\times (R_1,R_2),
\end{eqnarray}
where
\begin{equation}\label{J22}
J_2(u_0)(r)=-\frac{\Phi[\mathcal{B}u_0](r)}{\Phi[\mathcal{C}u_0](r)}\exp\!
\bigg[\!\int_{\!{r}}^{\!{R_2}}\frac{\Phi[\mathcal{B}u_0](\xi)}{\Phi[
\mathcal{C}u_0](\xi)}d\xi\bigg],\quad\forall\,r\in (R_1,R_2),
\end{equation}
and
\begin{equation}\label{q0}
q_0(t,r)\!:=J_2(u_0)(r)h_0(t)+N_3^0(u_0,u_1,g_1,f)(t,r),\;\quad
\forall\,(t,r)\in [0,T]\times (R_1,R_2).
\end{equation}
We have thus shown that the pair $(h,q)$ solves the fixed-point
system
$(\ref{hhh}),\,(\ref{q3})$.\\
We can summarize the result of this section in the following
equivalence theorem.
\begin{theorem}\label{3.2}
The pair  $(u,k)\in{\mathcal{U}}^{\,2,p}(T)\times
C^{\beta}\big([0,T];W^{1,p} (R_1,R_2)\big)$ is a solution to the
identification problem $\emph{\PIK}$, $\emph{H,K}\in
\{\emph{D,N}\}$, if and only if the triplet $(v,h,q)$ defined by
$(\ref{v})$ and $(\ref{q})$ belongs to $ {\mathcal{U}}_{\emph{H,K}}^{\,1,p}(T)
\times
C^{\beta} \big([0,T];\mathbb{R}\big)\times
C^{\beta}\big([0,T];L^{p}(R_1,R_2)\big)$ and solves problem
$(\ref{problem1})\!-\!(\ref{PSIV}), \,(\ref{hhh}),\,
(\ref{q3})$.
\end{theorem}
\section{An abstract formulation of problem\\ (\ref{problem1})--(\ref{PSIV}), (\ref{hhh}), (\ref{q3}).}
\setcounter{equation}{0}
Starting from the result of the previous section, we can
reformulate our identification problem in a Banach space
framework. Let $A:\mathcal{D}(A)\subset X \to X$ be a linear
closed operator satisfying the following assumptions:
\begin{itemize}
\item[(H1)]\emph{there exists $\zeta\in (\pi /2,\pi)$
such that the resolvent
set of $A$ contains $0$ and the open sector
${\Sigma}_{\zeta}=\{\mu\in \mathbb{C}:|\arg\mu|<\zeta\}$;}
\item[(H2)]\emph{there exists $M>0$ such that ${\|{(\mu I-A)}^{-1}\|}
_{\mathcal{L}(X)}\leqslant M|\mu{|}^{-1}$ for every $\mu\in
{\Sigma}_{\zeta}$.}
\item[(H3)]\emph{$X_1$ and $X_2$ are Banach spaces such that
$\mathcal{D}(A)=X_2\hookrightarrow X_1\hookrightarrow X $. Moreover,
$\mu\to {(\mu I-A)}^{-1}$ belongs to ${\cal L}(X;X_1)$ and satisfies the
estimate ${\|{(\mu I-A)}^{-1}\|}
_{\mathcal{L}(X;X_1)}\leqslant M|\mu{|}^{-1/2}$ for every $\mu\in
{\Sigma}_{\zeta}$.}
\end{itemize}
Here $\mathcal{L}(Z_1;Z_2)$ denotes, for any pair of Banach spaces
$Z_1$ and $Z_2$, the Banach space of all bounded linear operators
 from $Z_1$ into $Z_2$ equipped with the uniform-norm.
 In particular we set ${\cal L}(X)=\mathcal{L}(X;X)$.\\
By virtue of assumptions (H1), (H2) we can define the analytic
semigroup $\{{\rm{e}}^{tA}{\}}_{t\geqslant 0}$ of bounded linear
operators in $\mathcal{L}(X)$
 generated by $A$. As is well-known, there exist positive constants
$\widetilde{c_{k}}(\zeta)\; (k \in\mathbb{N})$ such that
$$
\|A^k {\rm{e}}^{tA}\|_{\mathcal{L}(X)}\leqslant
\widetilde{c_{k}}(\zeta)Mt^{-k}\,, \;\quad\forall t \in
{\mathbb{R}}_{+},\;\, \forall k\in\mathbb{N}.$$
 After endowing
$\mathcal{D}(A)$ with the graph-norm,
we can define the following family of interpolation spaces
${\mathcal{D}}_{A}(\beta,p)$, $\beta\in (0,1)$, $p\in [1,+\infty]$, which are
intermediate between $\mathcal{D}(A)$ and $X$:
\begin{eqnarray}\label{interpol1}
{\mathcal{D}}_{A}(\beta,p)=
\Big\{x\in X: |x|_{{\mathcal {D}}_{A}(\beta,p)} < +\infty\Big\}, \qq \mbox{if } p\in [1,+\infty],
\end{eqnarray}
where
\begin{equation}
{|x|}_{{\mathcal{D}}_{A}(\beta,p)} = \left\{
\begin{array}{l}
\ds \Big(\int_0^{+\infty}\!t^{(1-\beta)p-1}\|A{\rm e}^{tA}x\|_X^p\,dt
\Big)^{\! 1/p},\q \mbox{if } p\in [1,+\infty),
\\[5mm]
\sup_{0<t\le 1}\big(t^{1-\beta}\|A{\rm{e}}^{tA}x\|_X\big),\q \hskip 1.2truecm \mbox{if } p=\infty.
\end{array} \right.
\end{equation}
\pn
They are well defined by virtue of assumption (H1). Moreover, we set
\begin{equation}\label{interpol2}
{\mathcal{D}}_{A}(1+\beta,p)\!=\! \{x\in\mathcal{D}(A):Ax\in
{\mathcal{D}}_{A}(\beta,p)\}\,.
\end{equation}
Consequently, ${\mathcal{D}}_{A}(n+\beta,p)$, $n\in\mathbb{N},
\beta\in (0,1)$, $p\in [1,+\infty]$, turns out to be a Banach space when equipped
with the norm
\begin{equation}
{\|x\|}_{{\mathcal{D}}_{A}(n+\beta,p)}\!=\!
\sum_{j=0}^{n}{\|A^{j}x\|}_{X}+{|A^{n}x|}_{{\mathcal{D}}_{A}(\beta,p)}\,.
\end{equation}
In order to reformulate in an abstract form our identification
problem
$(\ref{problem1})\!-\!(\ref{PSIV})$,\,$(\ref{hhh})$,\,
$(\ref{q3})$
 we need the
following assumptions involving spaces, operators and data:
\begin{alignat}{13}
&(\textrm{H}4)\;\emph{$Y$ and $Y_1$ are Banach spaces such that
$Y_1\hookrightarrow Y$;}&
\nonumber\\[2mm]
&(\textrm{H}5)\;\emph{$B:\mathcal{D}(B)\subset X\rightarrow X$ is a linear
closed operator such that $X_2\subset \mathcal{D}(B)$;}&
\nonumber\\[2mm]
&(\textrm{H}6)\;\emph{$C:\mathcal{D}(C):=X_1\subset X\rightarrow
X$ is a linear closed operator;}&
\nonumber\\[2mm]
&(\textrm{H}7)\;\emph{$E\in\mathcal{L}(Y;Y_1)$, $\Phi\in\mathcal{L}(X;Y)$,
 $\Psi\in {X}^{\ast}$, ${\Phi}_1\in\mathcal{L}(X_1;Y)$,
${\Psi}_{1}\in
X_1^{\ast}$;}&
\nonumber\\[2mm]
&(\textrm{H}8)\;\emph{$\mathcal{M}$ is a continuous bilinear operator
from $Y\times {\wtil X}_1$ to $X$ and from $Y_1\times X$ to $X$,}
& \nonumber\\ &\qq\; \emph {where $X_1\hookto {\wtil X}_1$;}&
\nonumber\\[2mm]
&(\textrm{H}9)\;\emph{$J_1:X_2\rightarrow\mathbb{R}$,
$J_2:X_2\rightarrow Y$, $J_3:X_2\rightarrow\mathcal{L}(Y)$\,
are three prescribed (non-linear)}&
\nonumber\\
&\qq\;\emph{operators}\,;&\nonumber\\[1mm]
&(\textrm{H}10)\;\emph{$u_0, v_0\in X_2$,\, $Cu_0\in X_1$,\;
$J_1(u_0)\neq 0$, $Bu_0\in \mathcal{D}_A(\delta,+\infty)$, $\delta\in (\b,1/2)$ ;}&
\nonumber\\[2mm]
&(\textrm{H}11)\;\emph{$q_0\in C^{\beta}([0,T];Y)$,
$h_0\in  C^{\beta}([0,T])$\,;}
&\nonumber\\[2mm]
&(\textrm{H}12)\;\emph{$z_0\in C^{\beta}([0,T];X)$,\;
$z_1\in C^{\beta}([0,T];{\wtil X}_1)$,\; $z_2\in C^{\beta}([0,T];X)$\,;}&\nonumber\\[2mm]
&(\textrm{H}13)\;\emph{$Av_0+\mathcal{M}(\wtil{q}_0,Cu_0)+
\wtil{h}_0Bu_0-\mathcal{M}(E\wtil{q}_0,Bu_0)+z_2(0,\cdot)
\in\mathcal{D}_A(\beta,+\infty)$\,;}&\nonumber
\end{alignat}
where $\wtil{h}_0$ and $\wtil{q}_0$ are defined in the following Remark \ref{rem4.2}.\\
We can now reformulate our problem: \emph{determine a function
$v\in C^1([0,T];X)\cap C([0,T];X_2)$ such that}
\begin{eqnarray}
\label{problem2} v'(t)\!\!\!&=&\!\!\!{[\l_0I+A]}v(t)+\int_0^t\!\!\!
h(t-s)[{B}v(s)+z_0(s)]ds-\!
\int_0^t \mathcal{M}\big(Eq(t-s),{B}v(s)+z_0(s)\big)ds\nonumber\\[2mm]
&& + \int_0^t \mathcal{M}\big(q(t-s),{C}v(s)+z_1(s)\big)ds
+\mathcal{M}\big(q(t),{C}u_0\big)+h(t)Bu_0\nonumber\\[2mm]
& &-\mathcal{M}\big(Eq(t),Bu_0\big)+z_2(t),
\hskip 5truecm \forall\;t\in[0,T],\\[2mm]
\label{v02} v(0)\!\!\!&=&\!\!\!v_0.
\end{eqnarray}
\begin{remark}\label{z0z1z2}
\emph{In the explicit case $(\ref{problem1}), (\ref{PSIV})$ we have
$A={\cal{A}}-\l_0I$, with a large enough positive $\l_0$, and the
functions
 $z_0, z_1, z_2$ defined by}
\begin{eqnarray}
\label{z1z2z3}&z_0=D_t\mathcal{B}u_1\;,\qquad
z_1=D_t\mathcal{C}u_1\;,\qquad
 z_2=D_t\mathcal{A}u_1-D_t^2u_1+D_tf,&
\end{eqnarray}
\emph{whereas $v_0, h_0, q_0$ are defined, respectively,
 via the formulae $(\ref{3.3})$, $(\ref{h0})$, $(\ref{q0})$.}
\end{remark}
 Let us now introduce the following unknown function $w$ related to $v$ by
\begin{equation}\label{w}
w=Av\qquad\iff\qquad v=A^{-1}w\,.
\end{equation}
Applying $A$ to the Volterra operator equation equivalent to
problem $(\ref{problem2}),\,(\ref{v02})$ and using $(\ref{w})$, we
can easily obtain the following equation for $w$:
\begin{eqnarray}
\label{problem3} w(t)\!\!\!& =&\!\!\! A{\rm{e}}^{tA}v_0
+\l_0\int_0^t {\rm{e}}^{(t-s)A}w(s)ds +A\int_0^t
{\rm{e}}^{(t-s)A}z_2(s)ds \nonumber\\[2mm]
 & &\!\!\!+\, A\int_0^t {\rm{e}}^{(t-s)A}ds
\int_0^s h(s-\sigma)[{B}A^{-1}w(\sigma)+z_0(\sigma)]d\sigma
\nonumber\\[2mm]
& &\!\!\!-\, A\int_0^t {\rm{e}}^{(t-s)A}ds\int_0^s
\mathcal{M}\big(Eq(s-\sigma),BA^{-1}w(\sigma)
+z_0(\sigma)\big)d\sigma
\nonumber\\[2mm]
& &\!\!\!+\, A\int_0^t {\rm{e}}^{(t-s)A}ds\int_0^s
\mathcal{M}\big(q(s-\sigma),
CA^{-1}w(\sigma)+z_1(\sigma)\big)d\sigma
\nonumber\\[2mm]
 & &\!\!\!+\, A\int_0^t h(s){\rm{e}}^{(t-s)A}Bu_0ds
-A\int_0^t {\rm{e}}^{(t-s)A}\mathcal{M}\big(Eq(s),{B}u_0\big)ds
\nonumber\\[2mm]
& &\!\!\!+\, A\!\int_0^t\!\!{\rm{e}}^{(t-s)A}\mathcal{M}\big(q(s),
{C}u_0\big)ds\,,
\hskip 3truecm\forall\;t\in [0,T].
\end{eqnarray}
Denoting by $\mathcal{K}$ the convolution
operator
\begin{equation}\label{K}
\mathcal{K}(f,g):\!=\int_0^t\!\!\mathcal{M}\big(f(t-s),g(s)\big)ds\,,
\end{equation}
which maps $ C^{\beta}([0,T];Y_1)\times C([0,T];X)$ and
$ C^{\beta}([0,T];Y)\times C([0,T];X_1)$  into $C^{\beta}([0,T];X)$ (cf. \cite{CL}, section 4),
we can rewrite equation $(\ref{problem3})$ in the more compact way
\begin{equation}\label{problem4}
w=w_0+R_1(w,h,q)+S_1(q)\,,
\end{equation}
where we have set
\begin{align}
&R_1(w,h,q):= \l_0({\rm{e}}^{tA}\ast w)
+A[h\ast {\rm{e}}^{tA}\ast({B}A^{-1}w+z_0)]-A[ {\rm{e}}^{tA}
\ast \mathcal{K}(Eq,{B}A^{-1}w+z_0)]&\nonumber\\[2mm]
\label{R1}&\hskip 2truecm
+A[{\rm{e}}^{tA}\ast\mathcal{K}(q,{C}A^{-1}w+z_1)]
+A[{\rm{e}}^{tA} \ast
hBu_0] -A[{\rm{e}}^{tA}\ast
\mathcal{M}(Eq,{B}u_0)],&
\end{align}
and
\begin{eqnarray}
\label{S1}
S_1(q)\!\!\!&:=&\!\!\! A[{\rm{e}}^{tA}\ast \mathcal{M}(q,{C}u_0)],\\[2mm]
\label{w0}w_0\!\!\!&:=&\!\!\!A{\rm{e}}^{tA}v_0+A({\rm{e}}^{tA}\ast z_2).
\end{eqnarray}
Hence, applying operator $A^{-1}$ to both hand sides of
$(\ref{problem4})$, from $(\ref{R1})-(\ref{w0})$ we get
\begin{eqnarray}
v\!\!\!\!&=&\!\!\!\!A^{-1}w_0+\l_0({\rm{e}}^{tA}\ast v)
+{\rm{e}}^{tA}\ast h\ast ({B}A^{-1}w+z_0)
-{\rm{e}}^{tA}\ast
\mathcal{K}(Eq,{B}A^{-1}w+z_0)\nonumber\\[2mm]
&&\!\!\!\! +\, {\rm{e}}^{tA} \ast hBu_0+{\rm{e}}^{tA}\ast
\mathcal{K}(q,{C}A^{-1}w+z_1)-{\rm{e}}^{tA}\ast\mathcal{M}(Eq,{B}u_0)
+{\rm{e}}^{tA}\ast\mathcal{M}(q,{C}u_0)\nonumber\\[2mm]
\label{U} &&\!\!\!=: A^{-1}w_0+U(w,h,q)\,.
\end{eqnarray}
Now we rewrite the fixed point system $(\ref{hhh})$,\,$(\ref{q3})$
in the abstract form
\begin{eqnarray}
h(t)=h_0(t)\!\!\!\!&-&\!\!\!\!{\big[J_1(u_0)\big]}^{-1}
\Big\{\!\Psi\big[\mathcal{M}
\big(J_3(u_0)\{\Phi[{N_{1}}(v,h,q)(t)]-{\Phi}_1[v(t)]\},Cu_0\big)\big]
\nonumber\\[1,7mm]
&-&\!\!\!\!\Psi\big[\mathcal{M}\big(E\big(J_3(u_0)
\{\Phi[{N_{1}}(v,h,q)(t)]-{\Phi}_1[v(t)]\}\big),{B}u_0\big)\big]
\nonumber\\[1,7mm]
&-&\!\!\!\!\Psi[{N_{1}}(v,h,q)(t)]+{\Psi}_1[v(t)]\Big\}\nonumber\\[1,7mm]
\label{ha}&=:&\!\!\!\!h_0(t)+N_3(v,h,q)(t)\,,\;\;\,\qquad\qquad\qquad
\forall\,t\in[0,T],\\[3mm]
\label{qa} q(t)=q_0(t)\!\!\!\!&+&\!\!\!\!J_2(u_0)N_3(v,h,q)(t)+J_3(u_0)\{\Phi[{N_{1}}(v,h,q)(t)]-{\Phi}_1[v(t)]\}\,,
\nonumber\\[1,7mm]
&&\hskip 8truecm \forall\,t\in[0,T],
\end{eqnarray}
where $h_0$ and $q_0$ are the elements appearing in $(\textrm{H}11)$,
while $($cf. $(\ref{N1}))$ operator $N_1$ is defined by
\begin{equation}\label{N12}
N_{1}(v,h,q)(t)=-h\ast ({B}v+z_0)(t)+\mathcal{K}(Eq,{B}v+z_0)(t)-
\mathcal{K}(q,{C}v+z_1)(t)\,.
\end{equation}
\begin{remark}\label{rem4.2}
\emph{Since $N_{1}(v,h,q)(0)=0$, from $(\ref{ha})$ and $(\ref{qa})$ we can easily compute the
initial values $\wtil{h}_0$ and $\wtil{q}_0$ (appearing in $(\textrm{H}13))$ of functions $h$ and $q$:
\begin{equation}\label{rem4.2.1}
\left\{\!\begin{array}{l}
\wtil{h}_0=h_0(0)+J_4(u_0,v_0)=h(0)\,,\\[3mm]
\wtil{q}_0=q_0(0)+J_2(u_0)J_4(u_0, v_0)-J_3(u_0)\Phi_1[v_0]=q(0)\,,
\end{array}\right.
\end{equation}
where $J_4(u_0,v_0)$ is defined by:
\begin{eqnarray}
J_4(u_0,v_0)\!\!\!&=&\!\!\!{\big[J_1(u_0)\big]}^{-1}\!
\Big\{\!\Psi\big[\mathcal{M}
\big(J_3(u_0){\Phi}_1[v_0],Cu_0\big)-\mathcal{M}\big(EJ_3(u_0)
{\Phi}_1[v_0],{B}u_0\big)\big]\!-\!\Psi_1[v_0]\Big\}.\nonumber
\end{eqnarray}}
\end{remark}
\begin{remark}\label{rem4.3}
\emph{In the explicit case we get the equations
\begin{equation}
\wtil{h}_0=k_0(R_2)\,,\,\q\wtil{q}_0(r)=k_0'(r)\,.
\end{equation}
where $k_0$ is defined in $(\ref{k01})$.}
\end{remark}
Introducing the operators
\begin{align}
&\widetilde{R}_2(v,h,q)\!:=-{\big[J_1(u_0)\big]}^{-1}
\Big\{\!\Psi\big[\mathcal{M}\big(J_3(u_0)\Phi[{N_{1}}(v,h,q)],Cu_0\big)\big]&
\nonumber\\[1mm]
\label{tildeR2}&\qq\;\,\,\qquad\q-\Psi\big[\mathcal{M}\big(E\big(J_3(u_0)
\Phi[{N_{1}}(v,h,q)],{B}u_0\big)\big]-\Psi[{N_{1}}(v,h,q)]\!\Big\}\,,&
\end{align}
\begin{align}
\label{tildeR3}&\widetilde{R}_3(v,h,q)\!:=J_2(u_0)
\widetilde{R}_2(v,h,q)+J_3(u_0)\Phi[{N_{1}}(v,h,q)]\,,&
\\[1,8mm]
&\widetilde{S_2}(v)\!:=\!{\big[J_1(u_0)\big]}^{-1}
\Big\{\!\Psi\big[\mathcal{M}\big(J_3(u_0){\Phi}_1[v],Cu_0\big)\big]-
\Psi\big[\mathcal{M}\big(E\big(J_3(u_0){\Phi}_1[v],Bu_0\big)\big]
-{\Psi}_1[v]\!\Big\}\,,&\nonumber\\&&\label{tildeS2} \\
\label{tildeS3}&\widetilde{S_3}(v)\!:=J_2(u_0)
\widetilde{S_2}(v)-J_3(u_0){\Phi}_1[v]\,,&
\end{align}
the fixed-point system for $h$ and $q$ becomes
\begin{eqnarray}
\label{ha1}h\!\!\!&=&\!\!\!h_0+\widetilde{R}_2(v,h,q)+\widetilde{S_2}(v)\,,
\\[1,6mm]
\label{qa1}q\!\!\!&=&\!\!\!q_0+\widetilde{R}_3(v,h,q)+\widetilde{S_3}(v)\,.
\end{eqnarray}
Therefore, denoting
\begin{eqnarray}
\label{R_2}{R_2}(w,h,q)\!\!\!\!&:=&\!\!\!\widetilde{R}_2(A^{-1}w,h,q)\,,\\[1,6mm]
\label{R_3}{R_3}(w,h,q)\!\!\!\!&:=&\!\!\!\widetilde{R}_3(A^{-1}w,h,q)\,,
\end{eqnarray}
 and keeping in mind definitions $(\ref{R1})- (\ref{w0})$,\,
$(\ref{tildeR2})-(\ref{tildeS3})$, thanks to
$(\ref{problem4})$, $(\ref{ha1})$, $(\ref{qa1})$ we can pose
the following problem related to a given triplet $(w_0,h_0,q_0)\in
C^{\beta}([0,T];$ \newline $X)\times C^{\beta}([0,T]; \mathbb{R})\times
C^{\beta}([0,T];Y)$: \emph{determine a solution $(w,h,q) \in
C^{\beta}([0,T];X)\times C^{\beta}([0,T];$ \newline $\mathbb{R})\times
C^{\beta}([0,T];Y)$ to the fixed-point system}
\begin{equation}\label{problem5}
\left\{\!\!\begin{array}{lll}
w\!\!\!\!&=&\!\!\!\!w_0+R_1(w,h,q)+S_1(q)\,,\\[1,8mm]
h\!\!\!\!&=&\!\!\!\!h_0+R_2(w,h,q)+\widetilde{S_2}(A^{-1}w)\,,\\[1,8mm]
q\!\!\!\!&=&\!\!\!\!q_0+R_3(w,h,q)+\widetilde{S_3}(A^{-1}w)\,.
\end{array}\right.
\end{equation}
By  virtue of $(\ref{w})$, $(\ref{U})$ and the linearity of
$\widetilde{S_2},\, \widetilde{S_3}$ it is immediate to check that
system $(\ref{problem5})$ is equivalent to the following one:
\begin{equation}\label{problem6}
\left\{\!\!\begin{array}{lll}
w\!\!\!\!&=&\!\!\!\!w_0+R_1(w,h,q)+S_1(q)\,,\\[1,8mm]
h\!\!\!\!&=&\!\!\!\!\big[h_0+\widetilde{S_2}(A^{-1}w_0)\big]
+\big[R_2(w,h,q)+\widetilde{S_2}(U(w,hq))\big]
\\[1,4mm]
&=&\!\!\!\!\!:\overline{h}_0+R_5(w,h,q)\,,\\[1,8mm]
q\!\!\!\!&=&\!\!\!\!\big[q_0+\widetilde{S_3}(A^{-1}w_0)\big]
+\big[R_3(w,h,q)+\widetilde{S_3}(U(w,hq))\big]
\\[1,4mm]
&=&\!\!\!\!\!:\overline{q}_0+R_6(w,h,q)\,.
\end{array}\right.
\end{equation}
Hence, replacing $q$ in $S_1(q)$ with
$\overline{q}_0+R_6(w,h,q)$, and taking advantage of the
linearity of operator $S_1$, we deduce that the fixed-point
system $(\ref{problem6})$ is equivalent to the next one:
\begin{equation}\label{problem7}
\left\{\!\!\begin{array}{lll}
w\!\!\!\!&=&\!\!\!\!w_0+S_1(\overline{q}_0)+R_1(w,h,q)+S_1(R_6(w,h,q))\\[1,4mm]
&=&\!\!\!\!\!:\overline{w}_0+R_4(w,h,q)\,,\\[1,8mm]
h\!\!\!\!&=&\!\!\!\!\overline{h}_0+R_5(w,h,q)\,,\\[1,8mm]
q\!\!\!\!&=&\!\!\!\!\overline{q}_0+R_6(w,h,q)\,.\\
\end{array}\right.
\end{equation}
Finally, from $(\ref{w0})$ and the definitions of
$\overline{h}_0,\,\overline{q}_0$ in $(\ref{problem6})$ and of
$\overline{w}_0$ in $(\ref{problem7})$ we derive the following
representation of
$\overline{w}_0,\,\overline{h}_0,\,\overline{q}_0$ in terms of
$v_0,\,h_0,\,q_0,\,z_2$:
\begin{equation}\label{problem8}
\left\{\!\!\begin{array}{lll}
\overline{w}_0\!\!\!\!&=&\!\!\!\!A{\rm{e}}^{tA}v_0+A({\rm{e}}^{tA}\ast
z_2)+S_1(\overline{q}_0)\,,
\\[1,8mm]
\overline{h}_0\!\!\!\!&=&\!\!\!\!{h_0}
+\widetilde{S_2}({\rm{e}}^{tA}v_0+{\rm{e}}^{tA}\ast z_2)\,,
\\[1,8mm]
\overline{q}_0\!\!\!\!&=&\!\!\!\!{q_0}+
\widetilde{S_3}({\rm{e}}^{tA}v_0+{\rm{e}}^{tA}\ast z_2)\,.\\
\end{array}\right.
\end{equation}
Since the present situation is analogous to the one in \cite{CL}, we can
follow the same procedure used there (cf. sections 5 and 6) to get the
following local in time existence and uniqueness theorem.
\begin{theorem}\label{C-L}
Under assumptions $(\emph{H}1)-(\emph{H}13)$ there exists $T^{\ast}\in
(0,T)$ such that for any $\tau\in (0,T^{\ast}]$ the fixed-point
system $(\ref{problem7})$ has a unique solution $ (w,h,q)\in
C^{\beta}([0,\tau];X)\times C^{\beta}([0,\tau];\mathbb{R}) \times
C^{\beta}([0,\tau];Y)$.
\end{theorem}
An immediate consequence of Theorem $\ref{C-L}$ and of the
equivalence result proved in this section is the following
corollary.
\begin{corollary}
Under assumptions $(\emph{H}1)\!-\!(\emph{H}13)$ there exists $T^{\ast}\in
(0,T)$ such that for any $\tau\in (0,T^{\ast}]$ problem
$(\ref{problem2}), (\ref{v02}), (\ref{ha1}), (\ref{qa1})$ admits a
unique solution $(v,h,q)\in \big[C^{1+\beta}([0,\tau];X)\cap
C^{\b}([0,\tau];X_2)\big]\times C^{\beta}([0,\tau];\mathbb{R})\times
C^{\beta}([0,\tau];Y)$.
\end{corollary}
\section{Solving the identification problem (\ref{problem1})--(\ref{PSIV}), \newline(\ref{hhh}),
(\ref{q3}) and proving Theorem \ref{teoremaprincipe}}
\setcounter{equation}{0}
The basic result of this section is the following theorem.
\begin{theorem}\label{teoremaprincipeperv}
Let assumptions
$(\ref{unel}),\,(\ref{ipotesiaij})\!-\!(\ref{ipotesibijeci}),\,
(\ref{p3})\!-\!(\ref{primasuPsi})$ and condition $(\ref{HO})$ be
fulfilled. Moreover assume that the data enjoy the properties
$(\ref{richiestasuf})\!-\!(\ref{richiesteperg2})$,
inequalities $(\ref{J0}), (\ref{J1})$ and consistency conditions
$(\emph{C2,H,K})$, $(\ref{PHIV1})$, $(\ref{PSIV1})$.\\
 Then there exists
$T^{\ast}\!\in \!(0,T]$ such that the identification problem
$(\ref{problem1})\!-\!(\ref{PSIV})$, $(\ref{hhh})$,
$(\ref{q3})$  admits a unique
solution $(v,h,q)\in
{\mathcal{U}}_{\emph{H,K}}^{\,1,p}(T^{\ast})\times
C^{\beta}([0,T^{\ast}]; \mathbb{R})\times
C^{\beta}([0,T^{\ast}]; L^{p}(R_1,R_2))$  depending
continuously on the data with respect to the norms related to the
Banach spaces in $(\ref{richiestasuf})\!-\!(\ref{richiesteperg2})$.\\
In the case of the specific operators $\Phi$, $\Psi$ defined by
$(\ref{Phi1})$, $(\ref{Psi1})$ the previous result is still true
if we assume that $\lambda\in C^1(\partial\mbox{}B(0,R_2))$ and
$\psi\in C^1(\overline{\Omega})$ with $\psi\!=\!0$
 on the part
of $\partial\mbox{}\Om$ where the Dirichlet condition is
possibly prescribed.
\end{theorem}
\begin{proof}
For any $p\in (3,+\infty)$ let us choose the Banach spaces $X,{\wtil X}_1,X_1,X_2,Y,Y_1$
according to the rule
\begin{align}
& X=L^{p}(\Omega),\quad {\wtil X}_1=W^{1,p}({\Omega}),\quad X_1=W_{\textrm{H,K}}^{1,p}({\Omega}),
\quad X_2=W_{\textrm{H,K}}^{2,p}(\Omega),&
\\[1,7mm]
& Y=L^{p}(R_1,R_2),\quad Y_1=W^{1,p}(R_1,R_2),&
\end{align}
where the spaces $W_{\textrm{H,K}}^{1,p}({\Omega})$ are defined,
respectively, in $(\ref{WDD})-(\ref{WNN})$ with $\g=1/2$.\\
Of course, with this choice the operators
$\mathcal{B},\mathcal{C}$ defined by $(\ref{A})$ with
$\mathcal{D}(B)=X_2$, $\mathcal{D}(C)=X_1$, $Bu={\cal{B}}u$,
$Cu={\cal{C}}u$, satisfy assumptions $(\textrm{H}4)-(\textrm{H}6)$.\\
Let us define $A$ to be the second-order differential operator
${\cal{A}}-\l_0 I:{{\cal{D}}({A})}\subset L^p(\Om)\rightarrow L^p(\Om)$,
${\cal{A}}$ defined in $(\ref{A})$ and satisfying
$(\ref{unel})$, $(\ref{ipotesiaij})$
and $\l_0$ being any (fixed) positive constant.\\
To show that assumptions $(\textrm{H}1)-(\textrm{H}3)$ hold
we recall that $p\in (3, +\infty)$ and reason as in the
proof of theorem $7.3.6$ in \cite{PA}.
For this purpose we assume that $u\in W_{\textrm{H,K}}^{2,p}(\Om)$ is a
solution to the equation
\begin{equation}\label{5.10}
\l u-Au=f\,,\qq\q f\in L^p(\Om)\,.
\end{equation}
>From the identity
\begin{eqnarray}
\hskip -1truecm (\l+\l_0)\|u\|_{L^p(\Om)}^p \!\!\!&+&\!\!\!
\int_{\Om}\sum_{j,k=1}^3a_{j,k}(x)\big[(p-1)\g_j\g_k+
\delta_j\delta_k+ i(p-2)\g_k\delta_j\big]dx\nonumber\\[2mm]
\label{5.3}
\hskip -1truecm &=&\!\!\!\int_{\Om}f(x)|u(x)|^{p-2}\bar{u}(x)dx\,,\qq
f=\l u-Au\,,\qq \textrm{Re}\, \l \ge 0\,,
\end{eqnarray}
where $|u|^{(p-4)/2}\bar{u}D_ku = \g_k +i\delta_k$, $k=1,2,3$,
we easily derive the estimates
\begin{align}
\label{5.4}
&\hskip 0.5truecm \textrm{Re}(\l+\l_0)\,\|u\|_{L^p(\Om)}^p\le \|f\|_{L^p(\Om)}
\|u\|_{L^p(\Om)}^{p-1}\,, \qq\q \textrm{Re}\, \l \ge 0\,,&
\\[2mm]
\label{5.5}&\hskip 0.5truecm \int_{\Om}\sum_{k=1}^3\,\big(|\g_k|^2+|\delta_k|^2\big)dx\le
\frac{1}{\a_1}\|f\|_{L^p(\Om)}\|u\|_{L^p(\Om)}^{p-1}\,,&\\[2mm]
&\hskip 0.5truecm |\textrm{Im} (\l+\l_0)|\,\|u\|_{L^p(\Om)}^{p}\le
(p-2)\frac{\a_2}{2}\int_{\Om} \sum_{k=1}^3\,\big(|\g_k|^2+|\delta_k|^2\big)dx
+\|f\|_{L^p(\Om)}\|u\|_{L^p(\Om)}^{p-1}&\nonumber\\[2mm]
\label{5.6}
& \hskip 4,2truecm \le\big[(p-2)\frac{\a_2}{2}+1\big]
\|f\|_{L^p(\Om)}\|u\|_{L^p(\Om)}^{p-1}\,,\qq
\textrm{Re}\,\l\ge 0\,.&
\end{align}
>From $(\ref{5.4})$ and $(\ref{5.6})$ we deduce
\begin{equation}\label{5.7}
|\l+\l_0|\,\|u\|_{L^p(\Om)}\le \Big\{1+\big[(p-2)\frac{\a_2}
{2\a_1}+1\big]^2\Big\}^{1/2}\|f\|_{L^p(\Om)}\,,\qq\textrm{Re}\,
\l\ge 0\,,
\end{equation}
Therefore $(\l I-A)$ is injective and has a closed range
in $L^p(\Om)$ if $\textrm{Re}\,\l\ge 0$.\\ To show that
$(\l I-A)$ is also surjective if $\textrm{Re}\,\l\ge 0$, let $v\in L^{p'}(\Om)$,
$p'=p/(p-1)$, be
a function satisfying  $\int_{\Om} [\l u(x) -Au(x)]v(x)dx=0$ for any
$u\in W_{\textrm{H,K}}^{2,p}(\Om)$. From lemma 7.3.4 in \cite{PA}, which applies also to
our more general case, we deduce that $A$ is self-adjoint. Hence we get that $v\in
W_{\textrm{H,K}}^{2,p'}(\Om)$ and $\int_{\Om}u(x)[\bar{\l}v(x)-Av(x)]dx=0$
 for any $u\in W_{\textrm{H,K}}^{2,p}(\Om)$.
Since $W_{\textrm{H,K}}^{2,p}(\Om)$ is dense in $L^p(\Om)$
we deduce that $\bar{\l}v-Av=0$ in $L^{p'}(\Om)$,
$v\in W_{\textrm{H,K}}^{2,p'}(\Om)$. Then from the definition $A={\cal A} - \l_0I$
and the following inequality (cf. formula (7) in \cite{OK})
\[
{\rm Re}\, \langle {\cal A}v,|v|^{p'-2}{\ov v}\rangle \le 0,\qq
\forall v\in W^{2,p'}_{{\textrm {H,K}}}(\Om),
\]
where $p'\in (1,2)$, we easily conclude that $v=0$, i.e. the range of $(\l I-A)$ is
the entire space $L^p(\Om)$. Therefore $(\l I-A)$ is bijective for all $\l\in \mathbb{C}$
such that $\textrm{Re}\,\l\ge 0$ and as a consequence of
$(\ref{5.7})$  we have $\rho(A)\supset\{\l\in\mathbb{C}:
\textrm{Re}\,\l\ge 0\}$.\\
Finally, from proposition 2.1.11 in \cite{LU} and $(\ref{5.7})$
we deduce that $A$ is sectorial and its resolvent satisfies
the estimate
\begin{equation}\label{5.8}
\|(\l I- A)^{-1}\|_{{\mathcal{L}}(L^p(\Om))}\le
\frac{C_1}{|\l|}\,,\qq\forall\,\l\in{\Sigma}_{\zeta}\,,
\end{equation}
for some $\zeta\in (\pi/2,\pi)$. Hence $(\textrm{H}1)$
and  $(\textrm{H}2)$ hold.\\
Moreover, from $(\ref{5.7})$ with $\l = 0$ and theorem 3.1.1
in \cite{LU} we deduce the estimate
\begin{equation}\label{5.9}
\|u\|_{W^{2,p}(\Om)}\le C_2\|Au\|_{L^p(\Om)}\,,\qq
\forall\, u\in W_{\textrm{H,K}}^{2,p}(\Om)\,.
\end{equation}
Let now $u\in W_{\textrm{H,K}}^{2,p}(\Om)$ be a solution to the equation
(\ref{5.10}). Then, for any $\l\in {\Sigma}_{\zeta}$, we get
\begin{equation}\label{5.11}
\|u\|_{W^{2,p}(\Om)}\le C_2\|Au\|_{L^p(\Om)}\le C_2\big(
|\l|\|u\|_{L^p(\Om)}+\|f\|_{L^p(\Om)}\big)\le C_2(C_1+1)
\|f\|_{L^p(\Om)}\,.
\end{equation}
Finally, from the interpolation inequality
\begin{equation}\label{5.12}
\|u\|_{W^{1,p}_{\textrm{H,K}}(\Om)}\le C_3\|u\|_{W^{2,p}{\textrm{H,K}}(\Om)}^{1/2}
\|u\|_{L^{p}(\Om)}^{1/2}\le\frac{C_3[C_2(C_1+1)]^{1/2}}{|\l|^{1/2}}
\|f\|_{L^p(\Om)},\qq \forall u\in W^{2,p}_{\textrm{H,K}}(\Om)\,,
\end{equation}
we obtain that the resolvent $(\l I- A)^{-1}$ belongs to ${\cal{L}}(X;X_1)$
for any $\l\in {\Sigma}_{\zeta}$ and satisfies the estimate
\begin{equation}\label{5.13}
\|(\l I -A)^{-1}\|_{{\cal{L}}(X;X_1)}\le C_4|\l|^{-1/2}\,,
\qq\forall\,\l\in {\Sigma}_{\zeta}.
\end{equation}
Therefore $(\textrm{H}3)$ is satisfied, too.\\
Define now the operators $\Phi,{\Phi}_1,\Psi,{\Psi}_1$ respectively by
$(\ref{Phi1}),(\ref{phi1}),(\ref{Psi1}),(\ref{ps1})$ and
operators $E$ and
 $\mathcal{M}$ by
\begin{align}
\label{EE}&Eq(r)=\int_r^{R_2}\!\!q(\xi)d\xi\,,\qquad\forall\,r\in [R_1,R_2],
\\[1,8mm]
\label{MM}&\mathcal{M}(q,w)(x)=q(|x|)w(x)\,,\qquad\forall\;x\in\Omega\,.
\end{align}
Observe then that by virtue of H\"older's inequality we get
\begin{eqnarray}
\int_{R_1}^{R_2}{|Eq(r)|}^pdr\!\!\!&\leqslant&\!\!\!\int_{R_1}^{R_2}\!\bigg[
\int_{R_1}^{R_2}|q(\xi)|d\xi\bigg]^p dr\leqslant
\|q\|^p_{L^p(R_1,R_2)}
\int_{R_1}^{R_2}\! (R_2-R_1)^{p-1}dr\nonumber\\[1,7mm]
&=&\!\!{(R_2-R_1)}^p \|q\|^p_{L^p(R_1,R_2)},\qq \forall\,
q\in L^{p}(R_1,R_2).
\end{eqnarray}
Since $D_rEq(r)=-q(r)$, it follows that
$E\in\mathcal{L}\big(L^{^{p}}(R_1,R_2),W^{1,p}(R_1,R_2)\big)$. It
is an easy task to show that $\mathcal{M}$ is a continuous
bilinear operator from $L^{^{p}}(R_1,R_2)\times W^{1,p}({\Omega})$
to $L^{p}(\Omega)$ and from $W^{1,p}(R_1,R_2)\times L^{p}(\Omega)$
to $L^{p}(\Omega)$, since the Sobolev imbedding theorems with
$p\in (3,+\infty)$ imply $W^{1,p}(\Omega) \hookrightarrow C(\overline{\Omega})$
and $W^{1,p}(R_1,R_2)\hookrightarrow
C([R_1,R_2])$ . Hence assumptions $(\textrm{H}7)$ and
$(\textrm{H}8)$ are satisfied.\\
Moreover, the assumptions in $(\textrm{H}10)$, but $J_1(u_0)\neq 0$,
are satisfied according to $(\ref{richiestasuf})-(\ref{richiestaperAu0})$.
Then, if we define $J_1(u_0),\,J_2(u_0),\,J_3(u_0)$ according to
formulae $(\ref{J1}),(\ref{J22}), (\ref{N2})$, it immediately
follows that assumption $(\textrm{H}9)$ is satisfied as well as the
condition $J_1(u_0)\neq 0$ in $(\textrm{H}10)$ is.\\
Finally we estimate the vector $(v_0,z_0,z_1,z_2,h_0,q_0)$ in
terms of the
 data $(u_0,u_1,f,g_1,g_2)$. Definitions $(\ref{N10}),\,(\ref{N20}),\,
(\ref{N30}),\,(\ref{N0})$ imply
\begin{eqnarray}
N_1^0(u_1,g_1,f),\;N_3^0(u_0,u_1,g_1,f)\in
C^{\beta}([0,T];L^{^{p}}(R_1,R_2)),\nonumber\\[1,7mm]
N_2^0(u_1,g_2,f),\;N_0(u_0,u_1,g_1,g_2,f)\in
C^{\beta}([0,T]).\qq\qq\nonumber
\end{eqnarray}
Therefore from $(\ref{h0})$ and $(\ref{q0})$ we deduce
\begin{equation}
(h_0,q_0)\in  C^{\beta}([0,T])\times
C^{\beta}([0,T];L^{^{p}}(R_1,R_2)),
\end{equation}
whereas from $(\ref{z1z2z3})$ and hypotheses
 $(\ref{richiestasuf})\!-\!(\ref{richiestaperu1})$
 it follows
\begin{align}
&(z_0,z_1,z_2)\in C^{\beta}([0,T];L^{p}(\Omega))\times
C^{\beta}([0,T];W^{1,p}
(\Omega))\times C^{\beta}([0,T];L^{p}(\Omega)).&
\end{align}
Hence assumptions $(\textrm{H}11)\!-\!(\textrm{H}12)$ are also
satisfied.\\
To check condition $(\textrm{H}13)$ first we recall that in this case
the interpolation space ${\cal D}_A(\b,+\infty)$ coincides
with the Besov spaces $B_{\textrm{H,K}}^{2\b,p,\infty}(\Om)\!\equiv
\!{\big(L^p(\Omega), W_{\!\textrm{H,K}}^{2,p}(\Omega)\big)}_{\beta,\infty}$
(cf. \cite[\textrm{section 4.3.3}]{TR}). Moreover, we recall that
$B_{\!\textrm{H,K}}^{2\b,p,p}(\Omega)=W_{\!\textrm{H,K}}^{2\b,p}(\Omega)$.
Finally, we remind the basic inclusion (cf. \cite[section 4.6.1]{TR})
\begin{equation}\label{inclusion}
W^{s,p}(\Om)\hookto B^{s,p,\infty}(\Om)\,,\q\;\textrm{if}\;\,s\notin \mathbb{N}\,.
\end{equation}
Since our function $F$ defined in $(\ref{richiestaperA2u0})$ belongs
to $W_{\textrm{H,K}}^{2\b,p}(\Om)$, it is necessarily an element of
$B_{\textrm{H,K}}^{2\b,p\infty}(\Om)$. Therefore $(\textrm{H}13)$ is satisfied, too.
The proof is now complete.
\end{proof}
\textbf{Proof of Theorem \ref{teoremaprincipe}.} It easily follows from
Theorems \ref{3.2} and \ref{teoremaprincipeperv}.\ $\square$
\begin{remark}\label{su2.36}
\emph{We want here to give some insight into the somewhat involved
condition $(\ref{richiestaperA2u0})$. For this purpose we need to
assume that the coefficients $a_{i,j}\in C^{2+\a}(\ov{\Om})$,
$i,j=1,2,3$, for some $\a\in (2\b,1)$, $\b\in (0,1/2)$. First, we
give a sketch of how the membership of function $F$ in
$W_{\textrm{H,K}}^{2\b,p}(\Om)$ may be derived from $(\ref{k01})$
and the following stricter conditions on the linear operator
$\Phi_1$ appearing on $(\ref{quartasuPhi})$ and on the data
\begin{eqnarray}
&\Phi_1\in {\cal{L}}\big(W^{3,p}(\Om);W^{2,p}(\Om)\big)\,,&\\[2mm]
&{\cal{A}}u_0(\cdot)+f(0,\cdot)-D_tu_1(0,\cdot)\in W_{\rm{H,K}}^{2,p}(\Om)
\cap W^{3,p}(\Om)\,.&
\end{eqnarray}
We observe that in our application the function
appearing in $(\textrm{H}13)$ coincide, by virtue of formulae (3.20), (3.21),
with function $F$ defined in $(\ref{richiestaperA2u0})$.\\
To show that $F$ belongs to $W_{\textrm{H,K}}^{2\b,p}(\Om)$ we limit ourselves
to pointing out the following basic steps ensuring that function $k_0$
defined in (3.18) actually belongs to $C^{1+\a}([R_1,R_2])$:
\begin{itemize}
\item
[{\it{i)\ }}] for any $\rho\in C^{\a}(\ov\Om), \a\in (2\b,1),
w\in W^{2\b,p}(\Om)$, $\rho w \in W^{2\b,p}(\Om)\;$ and satisfies\\
\ the estimate $\|\rho w\|_{W^{2\b,p}(\Om)}\le C\|\rho\|_{C^{\a}(\ov\Om)}
\|w\|_{W^{2\b,p}(\Om)}$\,;
\item
[{\it{ii)\ }}] operator $\Phi$ maps $C^{\a}(\ov\Om)$
into $C^{\a}([R_1,R_2])$\,.
\end{itemize}
As for as the boundary conditions involved by assumption $(\textrm{H}13)$ are concerned,
we observe that they are missing when $(\textrm{H,K})=(\textrm{N,N})$, while in
the remaining case they are so complicated that we like better not to explicit them and we
limit to list them as
$$F\;\textrm{satisfies boundary conditions (H,K)}.$$
Of course, when needed, such conditions
can be explicitly computed in terms of the data and function $k_0$ defined in (3.18).}
\end{remark}
\section{The two-dimensional case}
\setcounter{equation}{0} In this section we deal with the planar
identification problem $\PIK$  related to the annulus
$\Omega=\{x\in\mathbb{R}^2\!:R_1<|x|<R_2\}$, $0<R_1<R_2$.\\
Operators $\mathcal{A}$, $\mathcal{B}$, $\mathcal{C}$ are defined
by $(\ref{A})$ simply replacing the subscript $3$ with $2$:
\begin{eqnarray}
\mathcal{A}\!=\!\!\sum_{j=1}^{2}D_{x_j}\big(\sum_{k=1}^{2}a_{j,k}(x)D_{x_k}
\big)\,,\quad\,\mathcal{B}\!=\!\!\sum_{j=1}^{2}D_{x_j}\big(\sum_{k=1}^{2}
b_{j,k}(x)D_{x_k}\big)\,,\quad\,
\mathcal{C}\!=\!\!\sum_{j=1}^{2}c_{j}(x)D_{x_j}\,.
\end{eqnarray}
Moreover, we assume that there exist two positive constants
${\alpha}_1$ and ${\alpha}_2$ with ${\alpha}_1\leqslant{\alpha}_2$
such that
\begin{equation}\label{unel3}
{\alpha}_1|\xi{|}^2\leqslant
\sum_{i,j=1}^{2}a_{ij}(x){\xi}_i{\xi}_j
\leqslant{\alpha}_2|\xi{|}^2\,,\qquad\,\forall\,
(x,\xi)\in\Omega\times \mathbb{R}^2.
\end{equation}
Furthermore we assume that the coefficients of operators
$\mathcal{A},\, \mathcal{B},\,\mathcal{C}\,$ satisfy also the
following properties corresponding to (\ref{RAD}),
(\ref{ipotesiaij}), (\ref{ipotesibijeci}):
\begin{eqnarray}
\label{ipotesiaij1}  a_{i,j}\in W^{2,\infty}({\Omega}),\q
a_{i,j}=a_{j,i},
\q b_{i,j}\in W^{1,\infty}(\Om),\q c_{i}\in L^{\infty}(\Omega),\;\,
i,j=1,2,\\[2mm]
\label {RAD2} \label{RAD1} \sum_{j,k=1}^2\,
x_jx_ka_{j,k}(x)=|x|^2h(|x|),\qq \forall x\in {\ov \Om},\q\qq\qq\qq
\end{eqnarray}
for some $h\in C({\ov \Om})$.
\pn In the present case an
example of admissible linear operators $\Phi$ and $\Psi$ is now
the following:
\begin{eqnarray}
\label{Phi12}
\hskip 0,74truecm\Phi [\!\!\!\!\! &v&\!\!\!\!\!](r)\!:=
\int_{\!0}^{2\pi}\!\!\!\!\lambda(R_2x')v(rx')d\varphi\,,\qq\\[1,7mm]
\label{Psi12} \hskip 0,74truecm\Psi[\!\!\!\!\! &v&\!\!\!\!\!]\!:=
\int_{\!R_1}^{R_2} r
dr\int_{\!0}^{2\pi}\!\!\!\!\psi(rx')v(rx')\,d\varphi,
\end{eqnarray}
where $(x_1,x_2)=(r\cos\!\varphi,r\sin\!\varphi)$,
 $x'=(\cos\!\varphi,\sin\!\varphi)$.\\
>From $(\ref{D123})$ we obtain
\begin{equation}\label{6.7}
\left\{\!\! \begin{array}{lll} D_{x_1}\!\!\! & =
&\!\!\!\cos\!\varphi D_{r}-\displaystyle\frac{\sin\!\varphi}
{r}D_{\varphi}\,,\\[5mm]
D_{x_2}\!\!\!  & = &\!\!\!{\sin\!\varphi}D_{r}+\displaystyle\frac{
\cos\!\varphi}{r}D_{\varphi}\,.
\end{array}\right.
\end{equation}
Therefore, setting $a_{i,j}(r,\varphi)=a_{i,j}(r\cos{\!\varphi},r
\sin{\!\varphi})$, from $(\ref{6.7})$ we deduce
\begin{eqnarray}
\label{6.8}\sum_{k=1}^{2}a_{1,k}(x)D_{x_k}\!\!\!& = &\!\!\!
f_1(r,\varphi)D_{r}+
\frac{f_2(r,\varphi)}{r}D_{\varphi}\,,\\
 \label{6.9}\sum_{k=1}^{2} a_{2,k}(x)D_{x_k}\!\!\!& = &\!\!\!  g_1(r,\varphi)D_{r}+
\frac{g_2(r,\varphi)}{r}D_{\varphi}\,,
\end{eqnarray}
functions $f_j,\,g_j$, $j=1,2$, being defined by
\begin{equation}\label{fj2}
\left\{\!\!
\begin{array}{lll}
f_1(r,\varphi)\!\!\!&:=&\!\!\!
{\widetilde{a}}_{1,1}(r,\varphi)\!\cos\!\varphi+
{\widetilde{a}}_{1,2}(r,\varphi)\!\sin\!\varphi\,,\\[1,7mm]
f_2(r,\varphi)\!\!\!&:=&\!\!\!\widetilde{a}_{1,2}(r,\varphi)\!\cos\!\varphi-
{\widetilde{a}}_{1,1}(r,\varphi)\!\sin\!\varphi\,,
\end{array}\right.
\end{equation}
\vskip -0,1truecm
\begin{equation}\label{gj2}
\left\{\!\! \begin{array}{lll}
g_1(r,\varphi)\!\!\!&:=&\!\!\!{\widetilde{a}}_{2,1}(r,\varphi)\!\cos\!\varphi+
{\widetilde{a}}_{2,2}(r,\varphi)\!\sin\!\varphi\,,\\[1,7mm]
g_2(r,\varphi)\!\!\!&:=&\!\!\!\widetilde{a}_{2,2}(r,\varphi)\!\cos\!\varphi-
{\widetilde{a}}_{2,1}(r,\varphi)\!\sin\!\varphi\,.
\end{array}\right.
\end{equation}
Hence, from $(\ref{6.7})-(\ref{6.9})$ we get
\begin{eqnarray}
\label{D1D1'} D_{x_1}\Big(\sum_{k=1}^{2}\,a_{1,k}(x)D_{x_k}\Big)
\!\!\! &=&\!\!\! D_{r}\Big[f_1(r,\varphi)\!\cos\!\varphi D_{r}
+ \frac{ f_2(r,\varphi)\cos\!\varphi}{r}D_{\varphi}\Big]\nonumber\\[1,3mm]
& &-\frac{ \sin\!\varphi}{r } D_{\varphi}\Big[
f_1(r,\varphi)D_r+\frac{ f_2(r,\varphi)}{r}D_{\varphi}\Big],\\[2mm]
\label{D2D2'} D_{x_2}\Big(\sum_{k=1}^{2}\,a_{2,k}(x)D_{x_k}\Big)
\!\!\! &=&\!\!\! D_{r}\Big[g_1(r,\varphi)\sin\!\varphi D_{r}
+\frac{ g_2(r,\varphi)\sin\!\varphi}{r}D_{\varphi}\Big]\nonumber \\[1,3mm]
&&-\frac{ \cos\!\varphi}{r }
D_{\varphi}\Big[g_1(r,\varphi)D_r+\frac{
g_2(r,\varphi)}{r}D_{\varphi}\Big].
\end{eqnarray}
Defining the following functions
\begin{equation}\label{kj}
k_j(r,\varphi) :\, = f_j(r,\varphi)\cos\!\varphi +
g_j(r,\varphi)\sin\!\varphi\,, \qquad j=1,2,
\end{equation}
and using $(\ref{fj2})\!-\!(\ref{gj2})$ we can easily check that,
by virtue of $(\ref{RAD1})$, we have
\begin{equation}
\label{hhhh}
k_1(r,\varphi)\!=\!\widetilde{a}_{1,1}(r,\varphi){\cos}^{2}\varphi+2\,
\widetilde{a}_{1,2}(r,\varphi)\!\cos\!\varphi\sin\!\varphi
+\widetilde{a}_{2,2}(r,\varphi){\sin}^{2}\varphi =: h(r).
\end{equation}
Then, rearranging the terms on the right-hand sides of
 $(\ref{D1D1'}), (\ref{D2D2'})$ we  obtain the following polar
representation for the second order differential operator
$\mathcal{A}$:
\begin{eqnarray}
\label{tildeA2} \widetilde{\mathcal{A}}\!\!\! & = &\!\!\!
D_r\Big[k_1(r)D_r+\frac{k_2(r,\varphi)}{r}D_{\varphi}\Big] -
\frac{\sin\!\varphi}{r}D_{\varphi}\Big[{f_1(r,\varphi)}D_r
+\frac{f_2(r,\varphi,\theta)}{r}D_{\varphi}\Big]\nonumber\\[1,5mm]
&& +\frac{\cos\!\varphi}{r}D_{\varphi}\Big[{g_1(r,\varphi)}D_r
+\frac{g_2(r,\varphi)}{r}D_{\varphi}\Big]\,.
\end{eqnarray}
\begin{remark}\emph{Similarly to the three-dimensional case
a class of coefficients $a_{i,j}$ satisfying property
$(\ref{RAD2})$ is
\begin{equation}\label{condsuaij1}
\left\{\begin{array}{lll}
a_{1,1}(x)\!\!\!&=&\!\!\!a(|x|)+\displaystyle
\frac{x_2^2[c(x)-b(|x|)]}{|x|^2}
+\displaystyle\frac{x_1^2d(|x|)}{|x|^2},\\[5mm]
a_{2,2}(x)\!\!\!&=&\!\!\!a(|x|)+\displaystyle
\frac{x_1^2[c(x)-b(|x|)]}{|x|^2}
+\displaystyle\frac{x_2^2d(|x|)}{|x|^2},\\[5mm]
a_{1,2}(x)\!\!\!&=&\!\!\! a_{2,1}(x)=\displaystyle\frac{\,x_1x_2[
b(|x|)-c(x)+d(|x|)]}{|x|^2},
\end{array}\right.
\end{equation}
where $a,b,d\in C^{2,\infty}([R_1,R_2])$ and
$c\in W^{2,\infty}({\Om})$, $a$ and  $c$
being, respectively, {\it positive} and {\it non-negative}
 functions such that
\begin{equation}\label{M111}
a(r)-b^+(r)-d^-(r)>0\,,\qq\forall\,r\in[R_1,R_2]\,.
\end{equation}
This property ensure the uniform ellipticity of $\mathcal{A}$.}
\end{remark}
\pn Working in  Sobolev spaces related to $L^p(\Om)$ with
\begin{equation}\label{P2}
p\in (2,+\infty)
\end{equation}
we note that our requirements on  operators $\Phi$ and $\Psi$ and
the data are the same as in $(\ref{primasuPhiePsi})\!-\!
(\ref{richiesteperg2})$ whereas the Banach spaces ${\mathcal{U}}^{s,p}(T)$ and
${\mathcal{U}}_{\,\textrm{H,K}}^{\,s,p}(T)$\, are still defined by
$(\ref{Us})$.
\begin{theorem}\label{teoremaprincipe1}
Let assumptions $(\ref{unel3})-(\ref{RAD2})$, $(\ref{P2})$,
$(\ref{primasuPhiePsi})-(\ref{primasuPsi})$ be fulfilled.
Moreover assume that the data enjoy the properties
$(\ref{richiestasuf})\! -\!(\ref{richiesteperg2})$
and  satisfy
inequalities $(\ref{J0}), (\ref{J1})$.\\ Then there exists
$T^{\ast}\in (0,T]$ such that the identification problem
$\emph\PIK\, (\emph{H,K}\in\{\emph{D,N}\}) $,  admits a unique
solution $(u,k)\in{\mathcal{U}}_{\,\emph{H,K}}^{\,2,p}( T^{\ast})\times
C^{\beta}\big([0,T^{\ast}],W^{1,p}(R_1,R_2)\big)$ depending
continuously on the data with respect to the norms pointed out in
$(\ref{richiestasuf})\!-\!(\ref{richiesteperg2})$.\\
In the case of the specific operators $\Phi$, $\Psi$ defined as in
$(\ref{Phi12}),\,(\ref{Psi12})$ the previous result is still true
if we assume $\lambda\in C^1(\partial\mbox{}B(0,R_2))$ and
  $\psi\in C^1(\overline{\Omega})$ with ${\psi}\!=\!0$ on the
part of $\partial\mbox{}\Om$ where the Dirichlet condition is
possibly prescribed. \end{theorem}
\begin{lemma}\label{PHIPSI1}
 When $\Phi$ and $\Psi$ are  defined by $(\ref{Phi12})$ and
$(\ref{Psi12})$, respectively, conditions
$(\ref{primasuPhiePsi})\!-\!(\ref{primasuPsi})$ are
 satisfied under assumptions $(\ref{ipotesiaij1}),\,(\ref{RAD1})$ on the
coefficients $a_{i,j}\;(i,j=1,2)$ and the hypotheses that
$\lambda\in C^1(\partial\mbox{}B(0,R_2))$ and $\psi\in
C^1(\overline{\Omega})$ with ${\psi}\!=\!0$  on the part of
$\partial\mbox{}\Om$ where the Dirichlet condition is
possibly prescribed.
\end{lemma}
\begin{proof}
It is essentially the same as that of Lemma \ref{PHIPSI}.
Therefore, we leave it to the reader.
\end{proof}
\section{Solving system (\ref{problemk1}) and (\ref{problemk2})}
\setcounter{equation}{0}
We solve here the following integro-differential system introduced in
remark 1.1, where $n\!=\!2,3$:
\begin{align}
\label{8.1}& \int_0^t \big\{ D_r k(t-s,r)D_tD_r u(s,r) +
k(t-s,r)[D_tD_r^2u(s,r) + (n-1)r^{-1}D_tD_ru(s,r)]\big\}ds&
\nonumber\\[2mm]
 &\qq + D_r k(t,r)D_r u(0,r) +
k(t,r)[D_r^2u(0,r)+ (n-1)r^{-1}D_ru(0,r)] = D_t{\wtil f}(t,r)\,,&
\nonumber\\[2mm]
\hskip -1.3truecm &\hskip 9truecm\forall\,
(t,r)\in [0,T]\times [R_1,R_2],&
\\
\label{8.2} & \int_{R_1}^{R_2}
\l(r)k(t,r)dr=g(t),\hskip 3.8truecm \forall\, t\in [0,T].&
\end{align}
We assume that the data $(u,f,g)$ enjoy the following properties:
\begin{alignat}{4}
\label{8.2.1}&u\in W^{1,1}\big((0,T);W^{2,p}(R_1,R_2)\big)\,,&\\[1,5mm]
\label{8.2.2}&\wtil{f}\in C^1\big([0,T];L^p(R_1,R_2)\big)\,,&\\[1,5mm]
\label{8.2.3}&g\in C\big([0,T];\mathbb{R}\big)\,,\\[1,5mm]
\label{8.2.4}&\l\in L^{p'}(R_1,R_2)\,,&
\end{alignat}
where $p\in [1, \infty]$ and $p'$ denotes the conjugate
exponent of $p$.\\
Like in section 3 we introduce the new unknowns
\begin{equation}
\label{8.3} h(t)=k(t,R_1),\qquad q(t,r)=D_rk(t,r),\qquad\forall\,(t,r)\in[0,T]\times[R_1,R_2]
\end{equation}
and express $k$ in terms of $h$ and $q$:
\begin{equation}
\label{8.4} k(t,r)=h(t)+\int_{R_1}^{r}
q(t,\rho)d\rho,\qquad\forall\,(t,r)\in[0,T]\times[R_1,R_2].
\end{equation}
Changing the order of integration, from (\ref{8.2}) and (\ref{8.4}) we immediately
derive the equation
\begin{eqnarray}
\label{8.5} && h(t)\int_{R_1}^{R_2} \l(r)dr+\int_{R_1}^{R_2}
\l_1(r)q(t,r)dr=g(t), \qq \forall\, t\in [0,T],
\end{eqnarray}
where
\begin{eqnarray}\label{8.6}
&& \l_1(r)=\int_{r}^{R_2} \l(\rho)d\rho,\qq \forall\, r\in
[R_1,R_2].
\end{eqnarray}
Assuming
\begin{eqnarray}
\label{8.7}
 \kappa^{-1}:=\int_{R_1}^{R_2} \l(r)dr \neq 0,\qq
\end{eqnarray}
from (\ref{8.5}) we easily deduce
\begin{eqnarray}
\label{8.8} && h(t)=\kappa g(t) - \kappa\int_{R_1}^{R_2}
\l_1(\rho)q(t,\rho)d\rho, \qq \forall\, t\in [0,T].
\end{eqnarray}
Assume now
\begin{eqnarray}
\label{8.9} && |D_r u(0,r)|\ge m >0,\qq \forall\, r\in [R_1,R_2].
\end{eqnarray}
Then, owing to (\ref{8.4}) and (\ref{8.8}), system (\ref{8.1}),
(\ref{8.2}) is equivalent to the following Volterra integral
equation of the second kind:
\begin{eqnarray}
\label{8.10} q(t,r)\!\!\!&-&\!\!\!\!
\kappa\a(r)\!\int_{R_1}^{R_2}\! \l_1(\rho)q(t,\rho)d\rho \,+
\,\a(r)\!\int_{R_1}^r\!q(t,\rho)d\rho \,+ \int_0^t\!
\b(t-s,r)q(s,r)ds
\nonumber\\[2mm]
&-&\!\!\!\! \kappa\!\int_0^t \!\g(t-s,r)ds\!\int_{R_1}^{R_2}\!
\l_1(\rho)q(s,\rho)d\rho\,+ \int_0^t
\!\g(t-s,r)ds\!\int_{R_1}^{r}\!q(s,\rho)d\rho\,
= \wtil{f}_1(t,r),\nonumber\\[2mm]
&&\qq\qq\qq\qq\qq\;\qq\qq\qq\qq\; \forall\, (t,r)\in [0,T]\times
[R_1,R_2],
\end{eqnarray}
where
\begin{eqnarray}
\label{8.11}
&& \a(r)=\frac{D_r^2u(0,r)+ (n-1)r^{-1}D_ru(0,r)}{D_ru(0,r)},\\[2mm]
\label{8.12}
&& \b(t,r)=\frac{D_tD_r u(t,r)}{D_ru(0,r)},\\[2mm]
\label{8.13}
&& \g(t,r)=\frac{D_tD_r^2u(t,r)+ (n-1)r^{-1}D_tD_ru(t,r)}{D_ru(0,r)},
\\[1mm]
\label{8.14} && \wtil{f}_1(t,r)=\frac{D_t{\wtil
f}(t,r)}{D_ru(0,r)}-\kappa g(t)\a(r) - \kappa\int_0^t
\g(t-s,r)g(s)ds\,.\qq\q
\end{eqnarray}
>From assumption $(\ref{8.2.1})\!-\!(\ref{8.2.3})$ we easily deduce that
\begin{eqnarray}
\label{8.14.1}\a\in L^p(R_1,R_2),
&\b\in L^1\big((0,T); C([R_1,R_2])\big),&
 \g\in L^1\big((0,T); L^p(R_1,R_2)\big),\qq \\[1.5mm]
\label{8.14.2}&\wtil{f}_1\in C\big([0,T]; L^p(R_1,R_2)\big)\,.&
\end{eqnarray}
Moreover $\a, \b, \g, \wtil{f}_1$ satisfy the estimates
\begin{alignat}{5}
\label{8.14.3}&\|\a\|_{L^p(R_1,R_2)}\leqslant
\frac{1}{m}\Big(1+\frac{n-1}{R_1}\Big)
\|u\|_{W^{1,1}((0,T); W^{2,p}(R_1,R_2))}\,,& \\[1.5mm]
\label{8.14.4}&\|\b(t,\cdot)\|_{L^p(R_1,R_2)}\leqslant
\frac{1}{m}\|D_tu(t,\cdot)\|_{W^{2,p}(R_1,R_2)}\,,&\\[1.5mm]
\label{8.14.5}&\|\g(t,\cdot)\|_{L^p(R_1,R_2)}\leqslant \frac{1}{m}
\Big(1+\frac{n-1}{R_1}\Big)\|D_tu(t,\cdot)\|_{W^{2,p}(R_1,R_2)}\,,&
\\[1.5mm]
&\|\wtil{f}_1(t,\cdot)\|_{L^p(R_1,R_2)}\leqslant \frac{1}{m}
\|D_t\wtil{f}(t,\cdot)\|_{L^{p}(R_1,R_2)}+|\kappa|\|g\|_{C([0,T];
 \mathbb{R})}\|\a\|_{L^p(R_1,R_2)}\nonumber\\
\label{8.14.6}&\qq\qq\qq\qq\;\;+|\kappa|\|g\|_{C([0,T]; \mathbb{R})}
\int_0^T\!\|\g(s,\cdot)\|_{L^p(R_1,R_2)}ds\,. &\qq
\end{alignat}
Consider now the auxiliary integral equation
\begin{eqnarray}
\label{8.15} \hskip -1.3truecm && q(t,r) +
\a(r)\!\int_{R_1}^{r}\!q(t,\rho)d\rho\,=\,f(t,r) + \kappa\a(r)
\!\int_{R_1}^{R_2}\! \l_1(\rho)q(t,\rho)d\rho\,.
\end{eqnarray}
Hence, setting
\begin{eqnarray}
\label{8.16}
Q(t,r)\!\!\!&=&\!\!\!\int_{R_1}^{r}\!q(t,\rho)d\rho\,,\\[1,5mm]
\label{8.17}z(t,r)\!\!\!&=&\!\!\!f(t,r) + \kappa\a(r)
\!\int_{R_1}^{R_2}\! \l_1(\rho)q(t,\rho)d\rho\,,
\end{eqnarray}
it turns out that $(\ref{8.15})$ is equivalent to the following
 Cauchy problem
\begin{equation}\label{8.18}
\left\{\!\!\begin{array}{l}
D_rQ(t,r)+\a(r)Q(t,r)=z(t,r), \qq\forall\,r\in (R_1,R_2),\\[2mm]
Q(t,R_1)=0,
\end{array}\right.
\end{equation}
which has the solution
\begin{equation}\label{8.19}
Q(t,r)=\int_{R_1}^r\!z(t,\rho)\exp\Big(\int_{r}^{\rho}\!\a(s)ds\Big)
d\rho\,.\q
\end{equation}
Therefore, using $(\ref{8.16})$, $(\ref{8.17})$ and replacing the expression for $Q$ in
$(\ref{8.19})$ into $(\ref{8.15})$, we get
\begin{equation}\label{8.20}
q(t,r)=f_1(t,r)\,
+\,\kappa\a(r)\Big(\int_{R_1}^{R_2}\!\l_1(\rho)q(t,\rho)d\rho\Big)\exp\Big(\int_r^{R_1}\!\!\a(s)ds\Big)\,,
\end{equation}
where
\begin{equation}\label{8.21}
f_1(t,r)=f(t,r)-\a(r)\!\int_{R_1}^r\!f(t,\rho)
\exp\Big(\int_r^{\rho}\!\a(s)ds\Big)d\rho\,.
\end{equation}
Multiplying then each side in (\ref{8.20}) by $\l_1(r)$ and
integrating over $[R_1,R_2]$, we easily derive the equation
\begin{eqnarray}
 \Big[1-\kappa\int_{R_1}^{R_2}\! \l_1(\rho)\a(\rho)
\exp\Big(\int_{\rho}^{R_1}\!\!\a(s)ds\Big)d\rho\Big]
\int_{R_1}^{R_2}\!\l_1(r)q(t,r)dr\,
=\,\int_{R_1}^{R_2}\! \l_1(r)f_1(t,r)dr\,,&&\nonumber\\[2mm]
\label{8.22}\forall\, t\in [0,T].\qq&&
\end{eqnarray}
By an integration by parts, which makes use of condition $(\ref{8.6})$
 and  $(\ref{8.7})$ we easily deduce the equality
\begin{equation}\label{8.22.1}
1-\kappa\int_{R_1}^{R_2}\! \l_1(\rho)\a(\rho)
\exp\Big(\int_{\rho}^{R_1}\!\!\a(s)ds\Big)d\rho =
\kappa\int_{R_1}^{R_2}\!\l(\rho)
\exp\Big(\int_{\rho}^{R_1}\!\!\a(s)ds\Big)d\rho\,.
\end{equation}
Assume now
\begin{eqnarray}
\label{8.23}
{\kappa}_{_1}=\int_{R_1}^{R_2}\!\l(\rho)
\exp\Big(\int_{\rho}^{R_1}\!\!\a(s)ds\Big)d\rho \neq 0.
\end{eqnarray}
Then from $(\ref{8.22})-(\ref{8.23})$ we deduce
\begin{eqnarray}
\label{8.24} && \int_{R_1}^{R_2}\!
\l_1(\rho)q(t,\rho)d\rho\,=\,(\kappa\kappa_{_1})^{-1}\!
\int_{R_1}^{R_2}\! \l_1(\rho)f_1(t,\rho)d\rho\,,\qq \forall\, t\in
[0,T].
\end{eqnarray}
Therefore, by easy
computations, from $(\ref{8.20})$ and $(\ref{8.21})$
 it follows that the solution to (\ref{8.15}) is given by
\begin{eqnarray}\label{8.25}
\q q(t,r)=f(t,r)+\int_{R_1}^{R_2}\!G(r,\rho)f(t,\rho)d\rho\,,
\qq\forall\,t\in[0,T],
\end{eqnarray}
where the Green function $G$ is defined as follows
\vskip 0,05truecm
\begin{equation}\label{8.26}
G(r,\rho)=\left\{\!\!\begin{array}{l}
\a(r)\bigg[\!\exp\Big(\displaystyle\int_r^{R_1}
\!\!\a(s)ds\Big)\displaystyle\int_{\rho}^{R_2}
\!\frac{\l(\sigma)}{\kappa_{_1}}
\exp\Big(\displaystyle\int_{\sigma}^{\rho}\!\a(s)ds\Big)
d\sigma\,-\,\displaystyle\exp\Big(\int_r^{\rho}\!\a(s)ds\Big)\!\bigg],
\\[5mm]
\qq\qq\qq\qq\qq\qq\qq\qq\qq\qq\qq\qq\;\;\; R_1\leqslant\rho
<r\leqslant R_2,
\\[4mm]
\,\a(r)\exp\Big(\displaystyle\int_r^{R_1}
\!\a(s)ds\Big)\displaystyle\int_{\rho}^{R_2}
\!\frac{\l(\sigma)}{\kappa_{_1}}
\exp\Big(\displaystyle\int_{\sigma}^{\rho}\!\a(s)ds\Big)
d\sigma,\;\q R_1\leqslant r <\rho \leqslant R_2.
\end{array}\right.
\end{equation}
Consequently, (\ref{8.10}) turns out to be equivalent to
\begin{eqnarray}\label{8.27}
 q(t,r) =\wtil{f}_1(t,r)\!\!\!&+&\!\!\!\int_{R_1}^{R_2}\!\!
G(r,\rho)\wtil{f}_1(t,\rho)d\rho\, +Lq(t,r)+
\int_0^t\! ds\!\int_{R_1}^{R_2}\!
G_1(t-s,r,\rho)q(s,\rho)d\rho ,\nonumber\\
\end{eqnarray}
where, denoted with $\displaystyle{\chi_{_{[R_1,r]}}}$ the
characteristic function of the interval $[R_1,r]$,  the operator
 $L$ and the function $G_1$ are defined, respectively, by the
following formulae
\begin{eqnarray}
\hskip -1.3truecm
Lq(t,r)\!\!\!&=&\!\!\!-\int_0^t \b(t-s,r)q(s,r)ds\,,\\[2mm]
 G_1(t,r,\rho)\!\!\!&=&\!\!\!\kappa\g(t,r)\l_1(\rho)\,-\,\g(t,r)
\displaystyle{\chi_{_{[R_1,r]}}}(\rho)-\,\b(t,\rho)G(r,\rho)
\nonumber\\[2mm]
\hskip -1.3truecm \label{8.28}&&\!\!\!+\kappa\l_1(\rho)\!
\int_{R_1}^{R_2}\!G(r,\xi)\g(t,\xi)d\xi\,-
\int_{\rho}^{R_2}\!G(r,\xi)\g(t,\xi)d\xi\,.
\end{eqnarray}
Observe now that, according to $(\ref{8.26})$, $G$ satisfies the
inequality
\begin{equation}\label{8.29}
\big|G(r,\rho)\big|\leqslant C_1|\a(r)|\qq\forall\,r,\rho\in(R_1,R_2)
\end{equation}
where $C_1\!=\!C_1(p,\kappa_1,\|\a\|_{L^1(R_1,R_2)},
\|\l\|_{L^1(R_1,R_2)})>0$.\\
Likewise $G_1$ satisfies the inequality
\begin{eqnarray}\label{8.30}
\|G_1(t,r,\rho)\big|\!\!\!&\leqslant&\!\!\!
\Big[1+|\kappa|\|\l\|_{L^1(R_1,R_2)}\Big]\Big[|\g(t,r)|+
C_1|\a(r)|\|\g(t,\cdot)\|_{L^p(R_1,R_2)}(R_2-R_1)^{1/p'}\Big]
\nonumber\\[2mm]
& &\!\!\!+\,C_1|\a(r)|\|\b(t,\cdot)\|_{L^{\infty}(R_1,R_2)}\,.
\end{eqnarray}
>From $(\ref{8.30})$ we easily deduce that
\begin{equation}
G_1\in L^1\big((0,T); L^p\big((R_1,R_2); L^{p'}(R_1,R_2)\big)\big)
\end{equation}
and satisfy the estimate
\begin{alignat}{3}
&\|G_1(t,\cdot,\cdot)\|_{L^p((R_1,R_2);
L^{p'}(R_1,R_2))}
\leqslant C_1\|\a\|_{L^p(R_1,R_2)}
\|\b(t,\cdot)\|_{L^{\infty}(R_1,R_2)}(R_2-R_1)^{1/p'}&
\nonumber\\[2mm]
&\;+\Big[1+|\kappa|\|\l\|_{L^1(R_1,R_2)}\Big]
\Big[1+C_1\|\a\|_{L^p(R_1,R_2)}(R_2-R_1)^{1/p'}\Big]
\|\g(t,\cdot)\|_{L^p(R_1,R_2)}(R_2-R_1)^{1/p'}&\nonumber\\[2mm]
\label{8.31}&\;:=l(t)\,.&
\end{alignat}
According to properties $(\ref{8.14.1})$, $(\ref{8.14.2})$ we easily
deduce
\begin{equation}\label{8.32}
l\in L^1\big((0,T);\mathbb{R}\big).
\end{equation}
Observe also that function
\begin{equation}\label{8.33}
w(t,r)=\wtil{f}_1(t,r)+\int_{R_1}^{R_2}\!G(r,\rho)\wtil{f}_1(t,\rho)dr\q
\end{equation}
belongs to $C\big([0,T]; L^p(R_1,R_2)\big)$. Moreover $w$ satisfies the
inequalities
\begin{eqnarray}\label{8.34}
|w(t,r)|\leqslant |\wtil{f}_1(t,r)|+C_1|\a(r)|
\|\wtil{f}_1(t,\cdot)\|_{L^p(R_1,R_2)}{(R_2-R_1)}^{1/p'},\qq\qq
\nonumber\\[1,5mm]
\hskip 7,5truecm \forall\,(t,r)\in (0,T)\times(R_1,R_2),
\end{eqnarray}
and
\begin{equation}\label{8.35}
\|w(t,\cdot)\|_{L^p(R_1,R_2)}\leqslant
\Big[1+C_1|\a\|_{L^p(R_1,R_2)}(R_2-R_1)^{1/p'}\Big]
\|\wtil{f}_1(t,\cdot)\|_{L^p(R_1,R_2)},\q\;\;\forall\,t\in (0,T)\,.\q
\end{equation}
Introduce now the Banach space $X\!=\!L^p(R_1,R_2)$ endowed with
the usual norm. Then we can rewrite the integral equation $(\ref{8.27})$
as the fixed-point equation
\begin{equation}\label{8.36}
q=w+L_1q
\end{equation}
where, according to $(\ref{8.27})$, the operator $L_1$ is defined via the following formula:
\begin{equation}\label{8.36.1}
L_1q(t,r)=Lq(t,r)+\int_0^tds\int_{R_1}^{R_2}
G_1(t-s,r,\rho)q(s,\rho)d\rho\,.
\end{equation}
Observe now that $L_1$ maps $C([0,T];X)$ into itself and satisfies the
following inequalities
\begin{eqnarray}
\|L_1q(t)\|_X\!\!\!&\leqslant&\!\!\!\int_0^t\!
\Big[\|\b(t-s,\cdot)\|_{C([R_1,R_2])}+
\|G_1(t-s,\cdot,\cdot)\|_{L^p((R_1,R_2); L^{p'}(R_1,R_2))}\Big]
\|q(s)\|_Xds
\nonumber\\[2mm]
\label{8.37}&\leqslant&\!\!\!\int_0^t\varphi(t-s)\|q(s)\|_Xds,
\qq\q\forall\,t\in(0,T)\,,
\end{eqnarray}
where $($cf. $(\ref{8.31}))$
\begin{equation}\label{8.38}
\varphi(t)=\|\b(t,\cdot)\|_{C([R_1,R_2])}+l(t),\qq\forall\,t\in (0,T),
\end{equation}
belongs to $L^1(0,T)$ according to properties $(\ref{8.14.1})$
and $(\ref{8.32})$.\\
To derive the estimate relative to $G_1$ observe that, owing to the
integral version of Minkowski's inequality, we get
\begin{eqnarray}
\Big\|\int_0^t\!ds\int_{R_1}^{R_2}\!G_1(t-s,\cdot,\rho)q(s,\rho)
d\!\!\!\!\!&\rho&\!\!\!\!\!
\Big\|_X\leqslant\Big\|\int_0^t\!\|q(s)\|_X
\Big(\int_{R_1}^{R_2}\!|G_1(t-s,\cdot,\rho)|^{p'}d\rho\Big)^{1/p'}ds
\Big\|_X\nonumber\\[1,5mm]
\label{8.39}\leqslant
\int_0^T\|\!\!\!\!\!&q&\!\!\!\!\!(s)\|_{X}\|G_1(t-s,\cdot,\cdot)\|
_{L^p((R_1,R_2); L^{p'}(R_1,R_2))}ds\,.
\end{eqnarray}
We introduce now in $C([0,T];X)$ the following weighted norm
\begin{equation}\label{8.40}
\|f\|_{\s}=\sup_{t\in[0,T]}\e^{-\s t}\|f(t)\|_X
\end{equation}
which is equivalent to the usual one.\\
Now from $(\ref{8.37})$ rewritten in the equivalent form
\begin{equation}\label{8.42}
\e^{-\s t}\|L_1q(t)\|_X\leqslant
\int_0^t\e^{-\s(t-s)}\varphi(t-s)\e^{-\s s}\|q(s)\|_Xds,
\qq\forall\,t\in(0,T)\,,
\end{equation}
and Young's theorem on convolution we deduce the basic estimate
\begin{equation}\label{8.43}
\|L_1q\|_{\s}\leqslant\|q\|_X\int_0^T\!\e^{-\s t}\varphi(t)dt\,.
\end{equation}
Consequently, the norm of the linear operator $L_1$ does not exceed
$\int_0^T\!\e^{-\s t}\varphi(t)dt$, which tends to $0$ as
$\s\rightarrow +\infty$. Therefore, if we choose a large enough $\s$,
then $I-L_1$ is continuously invertible according to Neumann's
theorem.\\
We have thus proved the following theorem
\begin{theorem}\label{teorema3}
let $u,\wtil{f}, g, \l$ satisfy properties $(\ref{8.2.1})-
(\ref{8.2.4})$ and let assumptions $(\ref{8.7})$, $(\ref{8.9})$,
$(\ref{8.23})$ be fulfilled. Then problem $(\ref{8.1})$, $(\ref{8.2})$
admits a unique solution $k\in C\big([0,T];\newline
 W^{1,p}(R_1,R_2)\big)$
continuously depending on $(u, \wtil{f}, g)$ with respect to the
norms pointed out.
\end{theorem}


\begin{thebibliography}{9999}
\bibitem[AD]{AD} Adams R. A.: Sobolev Spaces, Academic Press,
New York-San Francisco-London 1975.
\bibitem[CL]{CL} Colombo F., Lorenzi A.: {\it An identification problem related to
parabolic integrodifferential equations with non commuting spatial operators},
J. Inverse Ill Posed Problems, 8 (2000), 505--540.
\bibitem[JJ]{JJ} Janno J.: {\it An inverse problem arising in compression of visco-elastic
medium,} Proc. Estonian Acad. Sci. Phys. Math. 49 (2000), 75-89.
\bibitem[JW]{JW} Janno J., v. Wolfersdorf L.: {\it An inverse problem for identification
of a time- and space-dependent memory kernel of a special kind in heat conduction},
Inv. Prob. {\bf 15} (1999), pp. 1455-1467.
\bibitem[LU]{LU} A.\,Lunardi: Analytic semigroups and optimal regularity in
parabolic problems, Birkh\"auser Verlag, Basel 1995.
\bibitem[OK]{OK} Okazawa N.: {\it Sectorialness of second order elliptic operators
in divergence form}, Proc. Amer. Math. Soc. 113 (1991), 701-706.
\bibitem[PA]{PA} A.\,Pazy: Semigroups of linear operators and applications
to partial differential equations, Applied mathematical sciences vol. 44,
Springer-Verlag, New York 1983.
\bibitem[TR]{TR} Triebel, H.: Interpolation Theory, Function Spaces,
Differential Operators, North Holland Publ. Co., Amsterdam - New York - Oxford 1978.
\end{thebibliography}
\end{document}